\documentclass[8pt]{article}

\usepackage{amsfonts,,mathrsfs}
\usepackage{latexsym,amssymb}
\usepackage{amsmath, amsbsy}
\usepackage{amsopn, amstext}
\usepackage{graphicx, color, epstopdf}
\usepackage{threeparttable}
\usepackage{hyperref}

\usepackage{array}
\usepackage{bm}
\usepackage{multirow}
\usepackage{booktabs}
\usepackage{diagbox}
\usepackage{subcaption}

\def\R{\mathbb{R}}

\numberwithin{equation}{section}
\numberwithin{figure}{section}
\numberwithin{table}{section}
\numberwithin{Lem}{section}
\numberwithin{Defi}{section}
\numberwithin{Theo}{section}
\numberwithin{Rem}{section}
 \numberwithin{Coro}{section}

\allowdisplaybreaks[4]

\newcommand{\dx}{\,\mathrm{d}x}
\renewcommand{\i}{\mathrm{i}}
\def\q{\mathbf{q}}
\def\w{\mathbf{w}}
\def\p{\mathbf{p}}
\def\dy{\,\mathrm{d}y}

\def\L{\mathcal{L}}

\title{On perfectly matched layers of nonlocal wave equations in unbounded multi-scale media \thanks{This work is supported in NSFC under grants No. 11771035, 12071401 and NSAF U1930402, Natural Science Foundation of Hunan Province No. 2019JJ50572, Natural Science Foundation of Hubei Province No. 2019CFA007 and Xiangtan University 2018ICIP01.}}

\author{Yu Du\thanks{Department of Mathematics, Xiangtan University, Hunan, 411105, China({\tt duyu@xtu.edu.cn})}
\and Jiwei Zhang\thanks{School of Mathematics and Statistics, and Hubei Key Laboratory of Computational Science, Wuhan University, Wuhan 430072, China.
({\tt jiweizhang@whu.edu.cn})}}

\date{}

\begin{document}

\maketitle

\begin{abstract}
A nonlocal perfectly matched layer (PML) is formulated for the nonlocal wave equation in the whole real axis and numerical discretization is designed for solving the reduced PML problem on a bounded domain. The nonlocal PML poses challenges not faced in PDEs. For example, there is no derivative in nonlocal models, which makes it impossible to replace derivates with complex ones. Here we provide a way of constructing the PML for nonlocal models, which decays the waves exponentially impinging in the layer and makes reflections at the truncated boundary very tiny. To numerically solve the nonlocal PML problem, we design the asymptotically compatible (AC) scheme for spatially nonlocal operator by combining Talbot's contour, and a Verlet-type scheme for time evolution. The accuracy and effectiveness of our approach are illustrated by various numerical examples.


\end{abstract}

\section{Introduction}\label{sec_intro}
Over the last decades nonlocal models have attracted much attention owing to its potentially promising application in various disciplines of science and engineering, such as the preridynamic (PD) theory of continuum mechanics, and the modeling of nonlocal diffusion process, see \cite{bobaru2010the,du2012analysis,silling2000reformulation,weckner2005the,zhou2010mathematical}. Peridynamics originally introduced in \cite{silling2000reformulation} is a nonlocal formulation of elastodynamics which can more easily incorporate discontinuities such as cracks and damage, and has been extended past its original formulation including micropolar, nanofiber networks and so on \cite{bobaru2007influence,foster2010viscoplasticity,gerstle2007micropolar,weckner2009green}. While most nonlocal models are formulated on bounded domains with volume constraints, there are indeed applications in which the simulation in an infinite medium may be useful, such as wave or crack propagation in the whole space.

In this paper, we consider constructing perfectly matched layers (PMLs) to numerically solve the following nonlocal wave equation
\begin{align}
& (\partial_t^2+\mathcal{L}) q(x,t) = f(x,t), \quad x\in\R,\ t>0,\label{eq:nonlocalwave}\\
& q(x,0) = \psi_0(x), \quad \partial_t q(x,0) = \psi_1(x),\quad x\in\R, \label{eq:nonlocalwavecon}
\end{align}
where $q(x,t)$ represents the displacement field, $\psi_0(x)$ and $\psi_1(x)$ are the initial values, $f(x,t)$ is the body force. The nonlocal operator $\L$ acting on $q$ is defined by
\begin{align} 
	\mathcal{L} q(x,t) = \int_{\R} \big( q(x,t)-q(y,t) \big) \gamma\Big(y-x,\frac{x+y}{2}\Big) \mathrm{d}y,	\label{eq:nonlocalOperator}
\end{align}
where the nonnegative kernel function $\gamma(\alpha,\beta)$ satisfies
\begin{align}
\gamma(-\alpha,\beta)=\gamma(\alpha,\beta),\ \forall \alpha,\beta\in\R,\quad \mathrm{and}\quad \gamma(\alpha,\beta)=0,\ \mbox{if}\ |\alpha|>\delta>0. \label{eq:symkernel}
\end{align}
We assume the initial values $\psi_k(x)\ (k=0,1)$ and the source $f(x,t)$ are compactly supported over a bounded domain $\Omega_f$ for all $t$.


The aim of this paper is to develop an efficient numerical scheme to compute the solution of problem \eqref{eq:nonlocalwave}-\eqref{eq:nonlocalwavecon} on the whole real axis.  We are facing two difficulties: 
\begin{itemize}
\item The definition domain is unbounded. This requires us to construct artificial/absorbing boundary conditions (ABCs) which artificially bounds the computational domain without changing the solution of a PDE or nonlocal model. Here we consider perfectly matched layers (PMLs) of nonlocal models to overcome the unboundedness of spatial domain;
\item The kernel in the proposed nonlocal PML equation is complex-valued and depends on the time $t$. As a result, the modified nonlocal operator in the PML equation is given by a convolution in time,  which differs from the original nonlocal operator \eqref{eq:nonlocalOperator}. In addition, the simulations are implemented in multi-scale media. These require us to develop an asymptotically compatible (AC) scheme which should be consistent with both its local limiting model (i.e., taking $\delta \to 0$) and the nonlocal model itself (i.e., taking $\delta =\mathcal{O}(1)$). 
\end{itemize}

To overcome the first difficulty of the unboundedness of definition domain, the accurate ABCs is a successful approach by absorbing any impinging waves on the artificial boundaries/layers. The great progress has been made for the construction of ABCs for various nonlocal models, see \cite{DuHanZhangZheng,DuZhangZheng,ZYZD16,ZhengHuDuZhang,zheng2020stability}. In this paper, we will apply the perfectly matched layer (PML) to confine a bounded domain of physical interest. 
The PML has two important properties: (i) waves in the PML regions decay exponentially and; (ii) the returning waves after one round trip through the absorbing layer are very tiny if the wave reflects off the truncated boundary. These properties make it useful to simulate wave propagations in various media and fields, e.g., \cite{ber94,becache2004perfectly,berenger1996three,bermudez2007an,cw,cl03,collino1998the,turkel1998absorbing}. While the PML has been well developed for local problems, there are few works on PMLs for nonlocal problems \cite{antoine2019towards,DuZhangNonlocal1,DuZhangNonlocal2,wildman2011,wildman2012a,wang2013matching,ji2020artificial}. The main reason is that,  due to the nonlocal horizon, the design of PMLs for nonlocal models poses challenges not faced in the PDEs setting. For example, when constructing local PMLs, one replaces derivatives with respect to real numbers by the corresponding complex derivatives. However, this process cannot be applied to the nonlocal operator which is in the form of integral. 

In this paper, we provide a way of constructing an efficient PML for nonlocal wave problem \eqref{eq:nonlocalwave}--\eqref{eq:nonlocalwavecon}. To do so, we first reformulate the wave equation into a nonlocal Helmholtz equation by using Laplace transform. The Laplace transform introduces a complex variable $s$. After that, we apply the PML modifications, recently developed in \cite{DuZhangNonlocal1,DuZhangNonlocal2} for nonlocal Helmholtz equations, to derive PMLs for the resulting nonlocal Helmholtz equation with $s$. In this situation, the kernel is still analytically continued into the complex coordinates and consequently, the modified equation has a complex-valued kernel depending on complex value $s$. Finally, we transform the modified nonlocal equation into its time-domain form by inverse Laplace transform. As a result, we obtain the nonlocal wave equation with PML modifications.

In term of the discretization of the nonlocal PML equation, asymptotic compatibility (AC) schemes, a concept developed in \cite{TianDu,TianDu2}, is needed to discretize the nonlocal operator \cite{Du2016handbook}. 
In this paper, the kernel is taken by the following heterogeneous diffusion coefficient
\begin{equation}\label{sm}
0< \sigma(x)=\frac{1}{2}\int_\R s^2\gamma(s,x)ds < \infty,
\end{equation}
which implies that the nonlocal model is in multi-scale media. Under the assumption \eqref{sm}, the nonlocal operator \eqref{eq:nonlocalOperator} will converge
to a local operator \cite{DuZhangZheng} in the form of 
\begin{equation} \label{LO}
\lim_{\delta\to 0^+} \mathcal{L} q(x)=-\partial_x\left[\sigma(x)\partial_xq(x)\right].
\end{equation}
As $\delta\rightarrow 0$, the solution of problem \eqref{eq:nonlocalwave} will converge to the solution of local wave equation
\begin{align}
 \partial_t^2 q(x,t)-\partial_x\left[\sigma(x)\partial_x  \right] q(x,t)= f(x,t), \quad x\in\R,\ t>0. 
\end{align}

The AC scheme can ensure that numerical solutions of nonlocal models converge to the correct local limiting solution, as both the mesh size $h$ and the nonlocal effect $\delta$ tend to zero. One can refer to \cite{du2019asymptotically,tian2017a,TianDu,tian2014asymptotically} for more details of AC schemes. Noting that our nonlocal PML problem involves a new complex-valued kernel arising from the inverse Laplace transform, we here present the analogous ideas given in \cite{DuZhangZheng} to discretize the one-dimensional nonlocal operator with the general complex and time-dependent kernels and complex functions. For practical multi-scale simulations, we apply Talbot's contour \cite{weideman2006optimizing} to the inverse Laplace transform and obtain its approximation consisting of several sub-kernels. For each sub-kernel, we employ an AC scheme, developed for complex functions in \cite{DuZhangNonlocal1,DuZhangNonlocal2}, to discretize the resulting nonlocal operator. After that, we introduce some new auxiliary functions and reformulate the semi-discrete problem into a second-order ODE system, which is finally solved by a Verlet-type scheme.

The outline of this paper is organized as follows. In section~\ref{sec:NPML}, we design the nonlocal PMLs and obtain a truncated nonlocal PML problem on a bounded domain. In section~\ref{sec_dis}, we first spatially discretize the truncated nonlocal PML problem into an ODE system with the variable $t$ and solve it by a Verlet-type scheme. In section~\ref{sec_num} we introduce the basic setting of parameters for the discretization, and present numerical examples to verify the efficiency of the nonlocal PMLs and the convergence order of our numerical scheme.

\section{Nonlocal Perfectly Matched Layers} \label{sec:NPML}

We now consider the construction of nonlocal PMLs by using complex-coordinate approach. The complex-coordinate approach is essentially based on analytic continuation of the wave equation into complex spatial coordinates where the fields are exponentially decaying. To do so, we assume the initial data functions and the kernel functions satisfy the following properties:
\begin{itemize}
\item[A1:] $\psi_1$ and $f$ are compactly supported into a finite interval $\mathcal{D}=(x_l,x_r)$, and $\psi_0$ is compactly supported into $(x_l+\delta,x_r-\delta)$; 
\item[A2:] $\gamma$ is compactly supported over a strip $[-\delta,\delta]\times\R$ with $\delta\leq x_r-x_l$;
\item[A3:] $\gamma$ is homogeneous in both $[x_r,+\infty)$ and $(-\infty,x_l]$, namely,
\begin{align}
\gamma(\alpha,\beta) =& \gamma_L(\alpha),\quad \beta\in (-\infty, x_l+\delta/2],\\
\gamma(\alpha,\beta) =& \gamma_R(\alpha),\quad \beta\in (x_r-\delta/2,+\infty].
\end{align}
\end{itemize}
In the sequel we take $\gamma_L=\gamma_R=\gamma_\infty$ for brevity. 

Performing the Laplace transform on \eqref{eq:nonlocalwave}, we have 
\begin{align}
s^2\hat q(x,s) +\L \hat q(x,s) = \hat f(x,s)+s\psi_0(x)+\psi_1(x), \quad x\in\R, \label{eq:Hel}
\end{align}
where $\hat{q}(x,s) = \mathscr{L}(q(x,t);s)$ with $\mathscr{L}$ representing the Laplace transform in time with $\Re \{s\} >0$.

Noting the nonlocal operator $\mathcal{L}$ is self-adjoint, we can rewrite \eqref{eq:Hel} into the weak form of 
\begin{align*}
\int_{\R} s^2\hat q(x,s) v(x)\dx -& \frac12\int_{\R}\int_{\R} \big[\hat q(x,s) - \hat q(y,s) \big]\big[ v(x)-v(y) \big]\\
& \gamma\Big(y-x,\frac{x+y}{2}\Big)\dx\dy = \int_\R \big(\hat f(x,s)+s\psi_0(x)+\psi_1(x)\big)v(x)\mathrm dx,\quad \forall v\in C_0^\infty(\R).
\end{align*}
The PML modifications can be viewed as a complex coordinate stretching of the original problem by constructing an analytic continuation to the complex plane \cite{DuZhangNonlocal1,DuZhangNonlocal2}. In this paper, we take 
\begin{align}
	\tilde x:=\int_0^x \alpha(\eta,s)\mathrm{d}\eta= \int_0^x \Big(1 + \frac{z}{s}\sigma(\eta)\Big) \mathrm{d}\eta,\qquad \tilde y:=\int_0^y \alpha(\eta,s)\mathrm{d}\eta= \int_0^y \Big(1 + \frac{z}{s}\sigma(\eta)\Big) \mathrm{d}\eta, \label{eq:complexStreching}
\end{align}
where the absorption function $\sigma(\eta)\leq1$ is positive in $\R\setminus\mathcal{D}$ and is zero in $\mathcal{D}$. The PML coefficient $z$ is a real or complex constant, such as $z= 10$ or $ z = 10+\i$. By replacing
\begin{equation*}
x\to\tilde x(x,s),\quad y\to \tilde y(y,s),\quad
\dx\to\frac{\partial \tilde x}{\partial x}\dx=\alpha(x,s)\dx,\quad
\dy\to\frac{\partial \tilde y}{\partial y}\dy=\alpha(y,s)\dy,
\end{equation*}
we can transform Eq. \eqref{eq:Hel} into the following nonlocal equation with PML modifications
\begin{align*}
\int_{\R} s^2\hat q(\tilde x,s) v(\tilde x)\dx &- \frac12\int_{\R}\int_{\R} \big[\hat q(\tilde x,s) - \hat q(\tilde y,s) \big]\big[ v(\tilde x)-v(\tilde y) \big]\gamma\Big(\tilde y-\tilde x,\frac{\tilde x+\tilde y}{2}\Big)\alpha(x,s)\alpha(y,s)\dx\dy\\
& = \int_\R\big(\hat f(\tilde x,s)+s\psi_0(\tilde x)+\psi_1(\tilde x)\big) v(\tilde x)\alpha(x,s)\mathrm dx,\quad \forall v\in C_0^\infty(\R),
\end{align*}
which implies that 
\begin{align}
s^2\alpha(x,s)\hat q(\tilde x,s) + \int_{\R} \big[\hat q(\tilde x,s) - \hat q(\tilde y,s) \big] \gamma\Big(\tilde y-\tilde x,\frac{\tilde x+\tilde y}{2}\Big)\alpha(x,s)\alpha(y,s)\mathrm dy \notag\\
= \hat f(x,s)+s\psi_0(x)+\psi_1(x). \label{eq_HelPML}
\end{align}
Noting that to derive the right hand side of the above equation,  we have used the facts that $\tilde x=x,\alpha(x,s)=1$ for $x\in\mathcal{D}$, the initial data $\psi_k(k=0,1)$ and the source function $f$ are compactly supported in $\mathcal{D}$. Thus, we continue the equation~\eqref{eq:Hel} into \eqref{eq_HelPML} in complex coordinates. One can see that the solutions $\hat q(\tilde x,s)$ will not change in the interior domain $\mathcal{D}$ and exponentially decay in the absorbing region $\sigma(x)>0$ by choosing an appropriate PML coefficient $z$. 

We now perform the inverse Laplace transform to turn the equation back into the time-domain form. To do so, we set  
\begin{align}
\tilde q(x,t) = \mathscr{L}_s^{-1}[\hat q(\tilde x,s)],\qquad \tilde \gamma(x,y,t)=\mathscr{L}_s^{-1}\Big[\frac{1}{s}\gamma\Big(\tilde y-\tilde x,\frac{\tilde x+\tilde y}{2}\Big)\alpha(x,s)\alpha(y,s)\Big]. \label{eq_invkernel}
\end{align}
Since $\tilde x=x$ for $x\in\mathcal{D}$, we have $\tilde q(x,t)=q(x,t)$ for $x\in\mathcal{D}$ and all time $t$, which implies that $\tilde q(x,0)=q(x,0)$ and $\partial_t\tilde q(x,t)|_{t=0}=\partial_t q(x,t)|_{t=0}$ for $x\in \mathcal{D}$. Therefore, we can naturally assume that $\tilde q(x,0)=\psi_0(x)$ and $\partial_t\tilde q(x,0)=\psi_1(x)$. Then, we have the following inverse Laplace transforms 
\begin{align}
\mathscr{L}_s^{-1}[ s^2\alpha(x,s)\hat q(\tilde x,s)-s\psi_0(x)-\psi_1(x) ] &= \mathscr{L}_s^{-1}[ (s^2+zs\sigma(x))\hat q(\tilde x,s) -s\psi_0(x)-\psi_1(x)]\notag\\
& = \partial_t^2 \tilde q(x,t) + z\sigma(x) \partial_t \tilde q(x,t), \label{eq_wavePMLp1}
\end{align}
and
\begin{align}
& \mathscr{L}_s^{-1}\Big[\big[\hat q(\tilde x,s) - \hat q(\tilde y,s) \big]\gamma\Big(\tilde y-\tilde x,\frac{\tilde x+\tilde y}{2}\Big)\alpha(x,s)\alpha(y,s)\Big] \notag\\
=& \mathscr{L}_s^{-1}\Big[\big[\big(s\hat q(\tilde x,s)-q(x,0)\big) -\big(s \hat q(\tilde y,s)-q(y,0)\big) \big] \cdot \frac{1}{s} \gamma\Big(\tilde y-\tilde x,\frac{\tilde x+\tilde y}{2}\Big)\alpha(x,s)\alpha(y,s)\Big] \notag\\
& + \mathscr{L}_s^{-1}\Big[\big[q(x,0)-q(y,0) \big] \cdot \frac{1}{s} \gamma\Big(\tilde y-\tilde x,\frac{\tilde x+\tilde y}{2}\Big) \alpha(x,s)\alpha(y,s)\Big] \notag\\
=& \big[ \partial_t \tilde q(x,t)- \partial_t \tilde q(y,t) \big] \ast \tilde \gamma(x,y,t) +\big[q(x,0)-q(y,0) \big] \tilde \gamma(x,y,t),\label{eq_wavePMLp2}
\end{align}
where $*$ indicates the convolution of two functions in time. Combining \eqref{eq_wavePMLp1} and \eqref{eq_wavePMLp2} with \eqref{eq_HelPML} yields the nonlocal wave equation with PML modifications as
\begin{align}
\big(\partial_t^2 +z\sigma(x)\partial_t \big) \tilde q(x,t) + &\int_\R \big[ \partial_t \tilde q(x,t)- \partial_t \tilde q(y,t) \big] \ast \tilde \gamma(x,y,t) \mathrm dy \notag\\
&= f(x,t) - \int_\R\big[\psi_0(x)-\psi_0(y) \big] \tilde \gamma(x,y,t)\mathrm{d}y,\quad x\in\R.
\end{align}

Noting $ \tilde x =x, \tilde y = y \; (\forall x,y\in\mathcal{D})$ and $\mathrm{supp}\;\psi_0(x)\subset(x_l+\delta,x_r-\delta)$ (see A1), we have $\tilde \gamma(x,y,t)=\gamma\Big( y- x,\frac{ x+ y}{2}\Big)$, which implies that for all $x$,
\begin{align*}
\int_\R\big[\psi_0(x)-\psi_0(y) \big] \tilde \gamma(x,y,t)\mathrm{d}y = \int_{\mathcal{D}}\big[\psi_0(x)-\psi_0(y) \big]  \gamma\Big( y- x,\frac{ x+ y}{2}\Big)\mathrm dy.
\end{align*}

We finally have the nonlocal PML wave equations
\begin{align} \label{PMLw}
\big(\partial_t^2 +z\sigma(x)\partial_t \big) \tilde q(x,t) + \mathcal{L}_{pml} \partial_t \tilde q(x,t) = f(x,t) - \int_{\mathcal{D}}\big[\psi_0(x)-\psi_0(y) \big]  \gamma\Big( y- x,\frac{ x+ y}{2}\Big)\mathrm dy,
\end{align}
where the nonlocal operator $\mathcal{L}_{pml}$ for the PML is given by
\begin{align}
 \mathcal{L}_{pml}  \partial_t \tilde q(x,t)=\int_\R \big[ \partial_t \tilde q(x,t)- \partial_t \tilde q(y,t) \big] \ast \tilde \gamma(x,y,t) \mathrm dy.
\end{align}
We point out that the nonlocal PML operator $\mathcal{L}_{pml}$ involves a convolution in time, which differs from the original nonlocal operator $\mathcal{L}$.

Noting the PML equation \eqref{PMLw} is still defined on the whole space, we need to truncate the computational region at some sufficiently large $x$ by putting homogeneous Dirichlet boundary conditions. To do so, we define the PML layer $\mathcal{D}_p=(x_l-d_p,x_l]\cup[x_r,x_r+d_p)$ with the thickness $d_p$ of the absorbing layer, and define the boundary layer $\mathcal{D}_b$ of width $\delta$ which surrounds $\mathcal{D}\cup\mathcal{D}_p$.

Thus, we derive the following truncated nonlocal wave problem with PML modifications:
\begin{align}
&\big(\partial_t^2 +z\sigma(x)\partial_t \big) \hat{\tilde q}(x,t) +  \mathcal{L}_{pml}  \partial_t \hat{\tilde q}(x,t)\notag\\
&\qquad\qquad\qquad= f(x,t) - \int_{\mathcal{D}}\big[\psi_0(x)-\psi_0(y) \big]  \gamma\Big( y- x,\frac{ x+ y}{2}\Big)\mathrm dy,\quad x\in\mathcal{D}\cup\mathcal{D}_p, \label{eq_truPML1}\\
& \hat{\tilde q}(x,0) = \psi_0(x),\quad \partial_t \hat{\tilde q}(x,0) = \psi_1(x),\quad x\in\mathcal{D}_b\cup\mathcal{D}_p\cup\mathcal{D},\label{eq_truPML2}\\
& \hat{\tilde q}(x,t)=0,\quad x\in \mathcal{D}_b,\ 0<t\leq T \label{eq_truPML3}.
\end{align}
The solution $\hat{\tilde q}(x,t)$ above will use to approximate the solution $q(x,t)$ of problem \eqref{eq:nonlocalwave}--\eqref{eq:nonlocalwavecon} in $\mathcal{D}$.

\section{Discretization of the truncated nonlocal wave problem}\label{sec_dis}
We now consider the discrete scheme of problem~\eqref{eq_truPML1}--\eqref{eq_truPML3} by using the AC scheme, developed in \cite{TianDu,du2019asymptotically,tian2014asymptotically,tian2017a}, to discretize the nonlocal PML operator, and using the Verlet-type scheme to solve the discrete ODE system obtained from the spatial discretization.

We first introduce a spatial uniform grid $\{x_i\}_{i\in\mathbb{Z}}$ with mesh size $h$. For simplicity, we take $x_0=x_l-d_p$ and $x_{N+1}=x_r+d_p$, and set 
\begin{align}
&\mathcal{I}_p = \{i|1\leq i\leq N, x_i\in{\mathcal{D}}_p\},\quad\quad \mathcal{I} = \{i|1\leq i\leq N, x_i\in{\mathcal{D}}\}.
\end{align}
\subsection{The approximation of nonlocal PML operator $\mathcal{L}_{pml}$}
Here we consider the numerical approximation of the complex-valued kernel $\tilde \gamma(x,y,t)$ \eqref{eq_invkernel} given by
\begin{align}
\tilde \gamma(x,y,t)=\mathscr{L}_s^{-1}[\mathcal{K}(x,y,s) ]=\frac{1}{2\pi\i}\int_{\Gamma} \mathcal{K}(x,y,s)e^{st} \mathrm ds,\label{eq_invkernel2}
\end{align}
where $\mathcal{K}(x,y,s)$ represents
\begin{align}
&\mathcal{K}(x,y,s) = \frac{1}{s}\gamma\Big(\tilde y-\tilde x,\frac{\tilde x+\tilde y}{2}\Big)\alpha(x,s)\alpha(y,s). \label{eq:mk}
\end{align}

Denote by $\Omega_\mathcal{K}$ the $s$-complex domain where $\mathcal{K}(x,y,s)$ viewed as a function of the variable $s$ is analytic for any given $x,y\in \mathcal{D}_p$. The notation $\Gamma$ denotes the Bromwich line $\Re(s)=\eta$ initially, where the parameter $\eta$ is taken large enough such that the complement of $\Omega_\mathcal{K}$ lies in the half-plane $\Re(s)<\eta$. A typical approach of numerically approximating the inverse Laplace transform is to deform the Bromwich line into Talbot's contour \cite{weideman2006optimizing}
\begin{align}
\Gamma:\ s(\theta) = \omega +\mu(\theta\cot\theta+\nu\i\theta),\quad -\pi\leq\theta\leq\pi,
\end{align}
where $\omega,\mu$, and $\nu$ are real parameters such that $\Gamma$ encloses $\mathbb{C}\setminus\Omega_\mathcal{K}$.

We define the grid 
\begin{align*}
\theta_j=-\pi+\frac{\pi}{m}(2j-1),\ j=1,2,\cdots,m,
\end{align*}
and approximate Eq.~\eqref{eq_invkernel2} by the trapezoidal rule
\begin{align}
\tilde \gamma(x,y,t)=\frac{1}{2\pi\i}\int_{-\pi}^\pi \mathcal{K}(x,y,s(\theta))e^{s(\theta)t} s'(\theta) \mathrm d\theta \approx \sum_{j=1}^m \varpi_j\mathcal{K}(x,y,\xi_j)e^{\xi_j t}, \label{eq_tabapp}
\end{align}
where $\xi_j=s(\theta_j)$ are the sampling points on $\Gamma$ and $\varpi_j=\frac{s'(\theta_j)}{m\i}$ are associated quadrature weights.

By \eqref{eq_tabapp}, we derive the approximation of the nonlocal PML operator $\mathcal{L}_{pml}$
\begin{align}
\mathcal{L}_{pml} \partial_t \hat{\tilde q}(x,t) \approx& \int_{\R} \big[ \partial_t \hat{\tilde q}(x,t)- \partial_t \hat{\tilde q}(y,t) \big]   \ast \sum_{j=1}^m \varpi_j  \mathcal{K}(x,y,\xi_j)e^{\xi_jt} \mathrm dy \notag\\
=& \sum_{j=1}^m \varpi_j \Big(\int_{\R} \big[ \partial_t \hat{\tilde q}(x,t)- \partial_t \hat{\tilde q}(y,t) \big] \mathcal{K}(x,y,\xi_j) \mathrm dy\Big) \ast e^{\xi_jt} \notag\\
=:&  \sum_{j=1}^m  \mathcal{L}_{pml}^j \partial_t  \hat{\tilde q}(x,t)  \ast e^{\xi_jt},
\end{align}
which yields the approximation of Eq.~\eqref{eq_truPML1}
\begin{align}
\big(\partial_t^2 +z\sigma(x)\partial_t \big) \hat{\tilde q}(x,t) +& \sum_{j=1}^m \mathcal{L}_{pml}^j \hat{\tilde q}(x,t)\ast e^{\xi_jt}\notag\\
 =& f(x,t) -\int_{\mathcal{D}}\big[\psi_0(x)-\psi_0(y) \big]  \gamma\Big( y- x,\frac{ x+ y}{2}\Big)\mathrm dy,\quad x\in\mathcal{D}.
\end{align}

\subsection{The spatial discretization and semi-discrete problem}
We now consider the AC scheme, originally given in \cite{DuZhangZheng} and further developed for complex functions in \cite{DuZhangNonlocal1,DuZhangNonlocal2}, to discretize the nonlocal operators $\mathcal{L}_{pml}^j$. The approximation of $\mathcal{L}_{pml}^j $ is given as 
\begin{align}
\mathcal{L}_{pml}^{j,h} \partial_t\hat{\tilde q}(x_n,t) =& \varpi_j\int_{\R} \sum_{k\neq n}\phi_{x_k}(y)\frac{ \partial_t \hat{\tilde q}(x_n,t)- \partial_t \hat{\tilde q}(x_k,t) }{x_k-x_n}(y-x_n) \notag\\
&\qquad \cdot\mathcal{K}(\frac{x_n+x_k}{2}-\frac{y-x_n}{2},\frac{x_n+x_k}{2}+\frac{y-x_n}{2},\xi_j) \mathrm dy\label{eq_disker}\\
=& \sum_{k} \tilde a_{n,k}^j \partial_t \hat{\tilde q}(x_k,t),\qquad \forall n\in\mathcal{I}_p\cup\mathcal{I},  \quad j=1,\cdots,m, \notag
\end{align}
where $\tilde a_{n,n}^j = -\sum_{ k\neq n} \tilde a_{n,k}^j $ with 
\begin{align*}
	\tilde a_{n,k}^j = - \frac{\varpi_j}{(k-n)h} \int_{\R}\phi_{x_k}(y) (y-x_n) \mathcal{K}(x_n+\frac{x_k-y}{2},\frac{x_k+y}{2},\xi_j)\dy,\ k\neq n.
\end{align*}

Using the spatial discretizations \eqref{eq_disker}, we derive the following semi-discrete problem:
\begin{align}
&\big(\partial_t^2 +z\sigma(x_n)\partial_t \big)\hat{\tilde q}_n(t) + \sum_{j=1}^m\sum_{k}\tilde a_{n,k}^j \partial_t \hat{\tilde q}_k(t) \ast e^{\xi_jt}\notag\\
&\quad= f(x_n,t)-\int_{\mathcal{D}}\big[\psi_0(x_n)-\psi_0(y) \big]  \gamma\Big( y- x_n,\frac{ x_n+ y}{2}\Big)\mathrm dy,\quad 1\leq n\leq N,\ 0<t, \label{eq_sdp1}\\
& \hat{\tilde q}_n(0) = \psi_0(x_n),\quad \partial_t \hat{\tilde q}_n(0) = \psi_1(x_n),\quad n\in\mathbb{Z},\label{eq_sdp2}\\
& \hat{\tilde q}_n(t) = 0,\quad n<1\ \mathrm{or}\ n>N,\ 0< t\leq T,\label{eq_sdp3}
\end{align}
where $\hat{\tilde q}_n(t)\approx\hat{\tilde q}(x_n,t)$.

In the remainder, we consider the convolutions of functions $\partial_t \hat{\tilde q}_k(t)$ and $e^{\xi_jt}$ over the range $[0,t]$ by introducing the auxiliary functions
\begin{align}
p_{k,j}(t) = \partial_t \hat{\tilde q}_k(t) \ast e^{\xi_jt},\quad k\in\mathcal{I}_p\cup\mathcal{I},\ j=1,\cdots,m.
\end{align}
It's clear that functions $p_{k,j}$ satisfy the following ODEs
\begin{align}
\partial_t p_{k,j}(t) = \xi_j p_{k,j}(t) + \partial_t \hat{\tilde q}_k(t) \label{eq_pjode}
\end{align}
with the initial conditions $p_{k,j}(0) = 0$.
Then we reformulate the the semi-discrete problem~\eqref{eq_sdp1}--\eqref{eq_sdp3} into the following ODEs
\begin{align}
&\big(\partial_t^2 +z\sigma(x_n)\partial_t \big)\hat{\tilde q}_n(t) + \sum_{j=1}^m\sum_{k}\tilde a_{n,k}^j p_{k,j}(t) \notag\\
&\qquad = f(x_n,t)-\int_{\mathcal{D}}\big[\psi_0(x_n)-\psi_0(y) \big]  \gamma\Big( y- x_n,\frac{ x_n+ y}{2}\Big)\mathrm dy,\quad 1\leq n\leq N,\ 0<t,\label{eq_disPML1}\\
& \partial_t p_{n,j}(t) = \xi_j p_{n,j}(t) + \partial_t \hat{\tilde q}_n(t),\quad 1\leq n\leq N,\ 0<t\leq T,\label{eq_disPML2}\\
& \hat{\tilde q}_n(0) = \psi_0(x_n),\quad \partial_t \hat{\tilde q}_n(0) = \psi_1(x_n),\quad p_{n,j}(0)=0,\quad n\in\mathbb{Z},\label{eq_disPML3}\\
& \hat{\tilde q}_n(t) = 0,\quad p_{n,j}(t)=0,\quad n<1\ \mathrm{or}\ n>N,\ 0< t\leq T.\label{eq_disPML4}
\end{align}

\subsection{The Verlet-type ODE solver}

We here introduce the Verlet-type algorithm to numerically solve the ODE system \eqref{eq_disPML1}--\eqref{eq_disPML4}. Denote by $D_\sigma$ the $N\times N$  the diagonal matrix with entries $\sigma(x_1),\sigma(x_2),\cdots,\sigma(x_N)$, and by $\tilde A_j\ (j=1,2,\cdots,m)$ the $N\times N$ matrices with entries $\tilde a_{n,k}^j\ (n,k=1,2,\cdots,N)$. The the ODE system \eqref{eq_disPML1}--\eqref{eq_disPML4} can be rewritten into the following form
\begin{align}
\w(t)-\q'(t) & = 0,\\ 
\mathbf{w}'(t) + zD_\sigma \mathbf{w}(t) + \sum_{j=1}^m \tilde A_j \mathbf{p}_j(t) &= \mathbf{f}(t),\\
\mathbf{p}_j'(t) - \xi_j \mathbf{p}_j(t) - \mathbf{w}(t) &= 0,\quad j=1,2,\cdots,m,
\end{align}
where $\mathbf{q}=(\hat{\tilde q}_1,\hat{\tilde q}_2,\cdots,\hat{\tilde q}_N)^T$, $\mathbf{p}_j=(p_{1,j},p_{2,j},\cdots,p_{N,j})^T$ and
\begin{align*}
\mathbf{f} = \begin{pmatrix}
f(x_1,t)-\int_{\mathcal{D}}\big[\psi_0(x_1)-\psi_0(y) \big]  \gamma\Big( y- x_1,\frac{ x_1+ y}{2}\Big)\mathrm dy\\
f(x_1,t)-\int_{\mathcal{D}}\big[\psi_0(x_2)-\psi_0(y) \big]  \gamma\Big( y- x_2,\frac{ x_2+ y}{2}\Big)\mathrm dy\\
\vdots\\
f(x_1,t)-\int_{\mathcal{D}}\big[\psi_0(x_N)-\psi_0(y) \big]  \gamma\Big( y- x_N,\frac{ x_N+ y}{2}\Big)\mathrm dy
\end{pmatrix}.
\end{align*}
Let $\tau$ be the temporal stepsize and $t_k=k\tau$ be the $k$-th time point. Denote by $
\w^{k+1/2} \approx \w(t_{k+1/2}),\; \q^k \approx \q(t_k),\; \p_j^k \approx \p(t_k)\; (k=0,1,\cdots).$ 
Let $\mathbf{f}^j=\mathbf{f}(t_j), \mathbf{\Psi}_0 =  (\psi_0(x_1),\psi_0(x_2),\cdots,\psi_0(x_N))^T$ and $\mathbf{\Psi}_1 =  (\psi_1(x_1),\psi_1(x_2),\cdots,\psi_1(x_N))^T$.
The initial values can be written as
\begin{align} \label{1s}
\q^0 &= \mathbf{\Psi}_0,\\
\p_j^0 &=0,\quad j=1,2,\cdots,m,\\
\w^{1/2} &= \mathbf{\Psi}_1 + \frac{\tau}{2}\Big[\mathbf{f}^0-\sum_{j=1}^m\tilde A_j \p_j^0-zD_\sigma\mathbf{\Psi}_1\Big].
\end{align}
For $k\geq0$, we apply the following second-order central difference to calculate $\q^{k+1}$ and $\p_j^{k+1}$ by 
\begin{align}
\w^{k+\frac12} - \frac{\q^{k+1}-\q^{k}}{\tau} & =0,\\
\frac{\p_j^{k+1}-\p_j^{k}}{\tau} - \xi_j\frac{\p_j^{k+1}+\p_j^{k}}{2} - \w^{k+\frac12} &=0.
\end{align}
After that, we still update $\w^{k+3/2}$ by using the above $\q^{k+1}$ and $\p_j^{k+1}$ via the second-order scheme: 
\begin{align}
& \frac{\w^{k+\frac32}-\w^{k+\frac12}}{\tau} + zD_\sigma \frac{\w^{k+\frac32}+\w^{k+\frac12}}{2} + \sum_{j=1}^m \tilde A_j \mathbf{p}_j^{k+1}= \mathbf{f}^{k+1}. \label{3s}
\end{align}

\section{Numerical examples}\label{sec_num}
In this section, three examples are provided to verify the effectiveness of our PML strategy,  the convergence and asymptotic compatibility of scheme \eqref{1s}--\eqref{3s}. Define the $L^2$-error at $t=t_k$ by 
\begin{align}
e_h &= \sqrt{\frac{1}{|\mathcal{I}|}\sum_{n\in\mathcal{I}} |\hat{\tilde q}^k(x_n)-q(x_n,t_k)|^2 },\end{align}
and the error to study the AC property, i.e., the so-called ``$\delta$-convergence'' in \cite{silling2000reformulation,DuZhangZheng,bobaru2009convergence}, by 
\begin{align}
e_\delta &= \sqrt{\frac{1}{|{\mathcal{I}\cup\mathcal{I}_p}|}\sum_{n\in{\mathcal{I}\cup\mathcal{I}_p}} |\hat{\tilde q}^k(x_n)-u(x_n,t_k)|^2 },
\end{align}
where $u(x,t)$ is the corresponding local PML soltion.

In the simulations, we choose the interior domain $\mathcal{D}$ such that $x_l=-l,\ x_r=l$ for some constant $l$ and set the PML absorbing function as the piecewise linear function
\begin{align}
\sigma(\eta) = 
\begin{cases}
0, & -l<\eta<l,\\
\frac{1}{d_{p}} (|\eta|-l), & l\leq|\eta|<d_{p}+l,\\
1, &  d_{p}+l\leq |\eta|.
\end{cases}
\end{align}

\noindent\textbf{Example 5.1.} In this example we take the source $f(x,t)\equiv0$ and the initial values as 
\begin{align*} 
\psi_0(x) = e^{-20(x-0.2)^2} + e^{-20(x+0.2)^2}, \quad \psi_1(x) = 100x^2e^{-20x^2},
\end{align*}
and consider the kernel
\begin{align}
\gamma(y-x,\frac{x+y}{2})=\frac{1}{\delta^3}\gamma_0(\frac{|x-y|}{\delta}),\label{eq_ker1}
\end{align}
where $\gamma_0(s) = \frac{1}{2c_0^3} e^{-\frac{|s|}{c_0}}$ and $c_0$ can be taken as a positive constant measuring the nonlocal horizon. 

In the simulations, we take $l=1.5$, $d_p=1$, the PML coefficient $z=20$, and set the parameters in the numerical implementation of Laplace transform as 
$\mu = \frac12 |z|, \nu = 1.$
We leave the discussion about the domain $\Omega_\mathcal{K}$ for the kernel~\eqref{eq_ker1} in the appendix~\ref{sec_Ap}.

\begin{figure}[htbp]
\centering
\begin{subfigure}{.32\textwidth}
 \centering
	\includegraphics[width=1\textwidth]{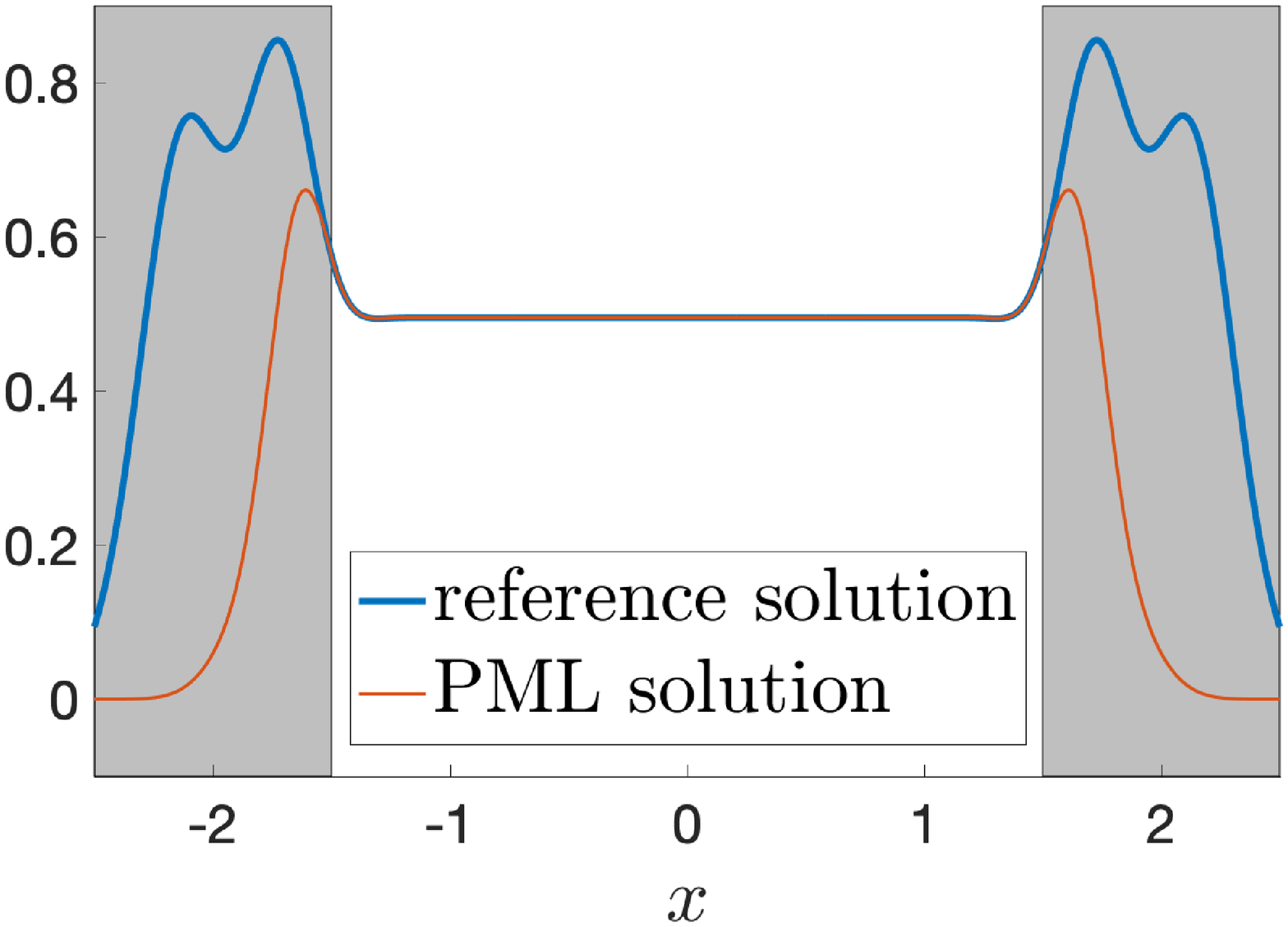}
 \caption{$\delta=0.5$}
 \label{fig:ex5p1solutionsT2a}
\end{subfigure}
\begin{subfigure}{.32\textwidth}
 \centering
	\includegraphics[width=1\textwidth]{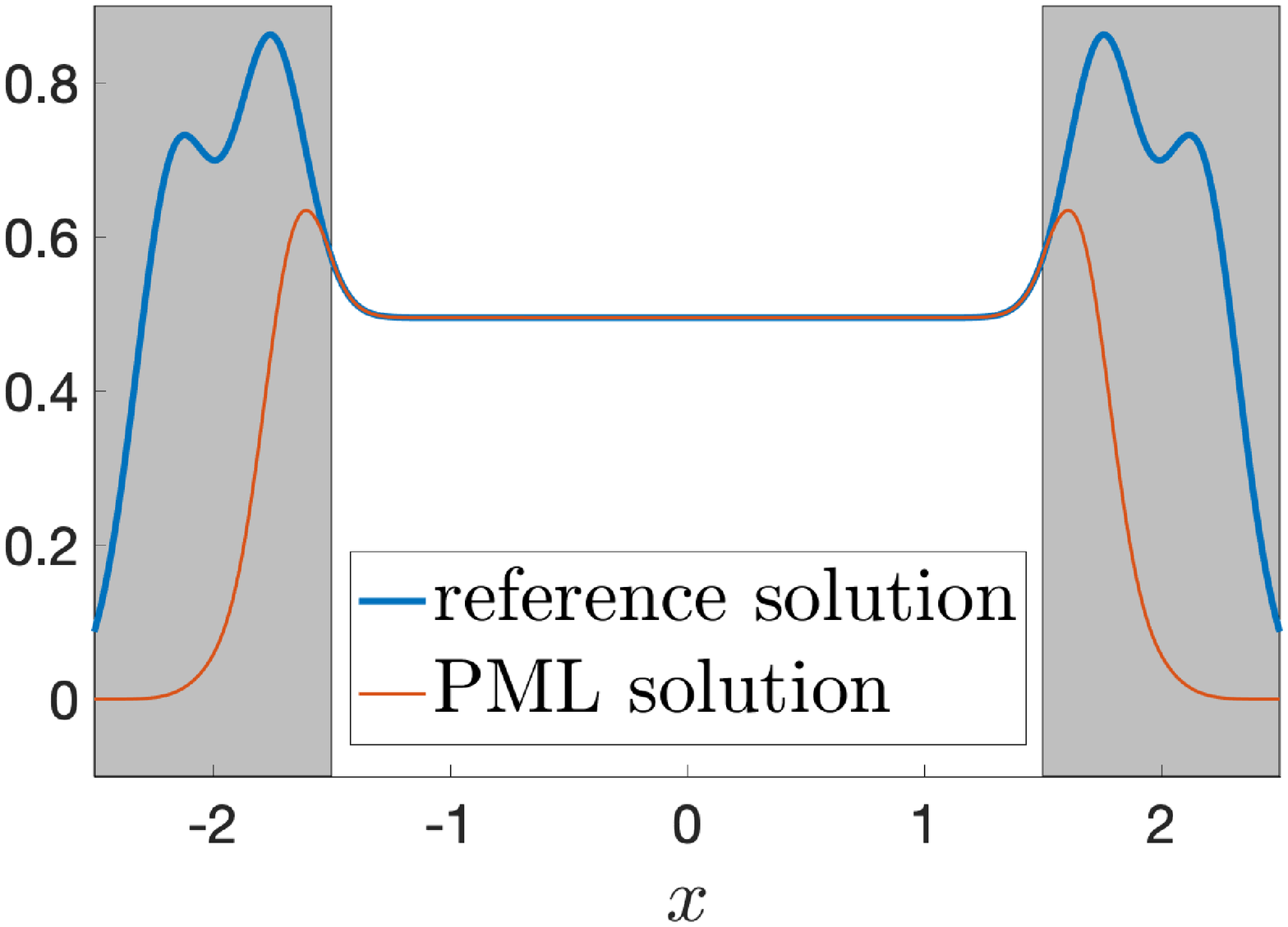}
 \caption{$\delta=0.2$}
 \label{fig:ex5p1solutionsT2b}
\end{subfigure}
\begin{subfigure}{.32\textwidth}
 \centering
	\includegraphics[width=1\textwidth]{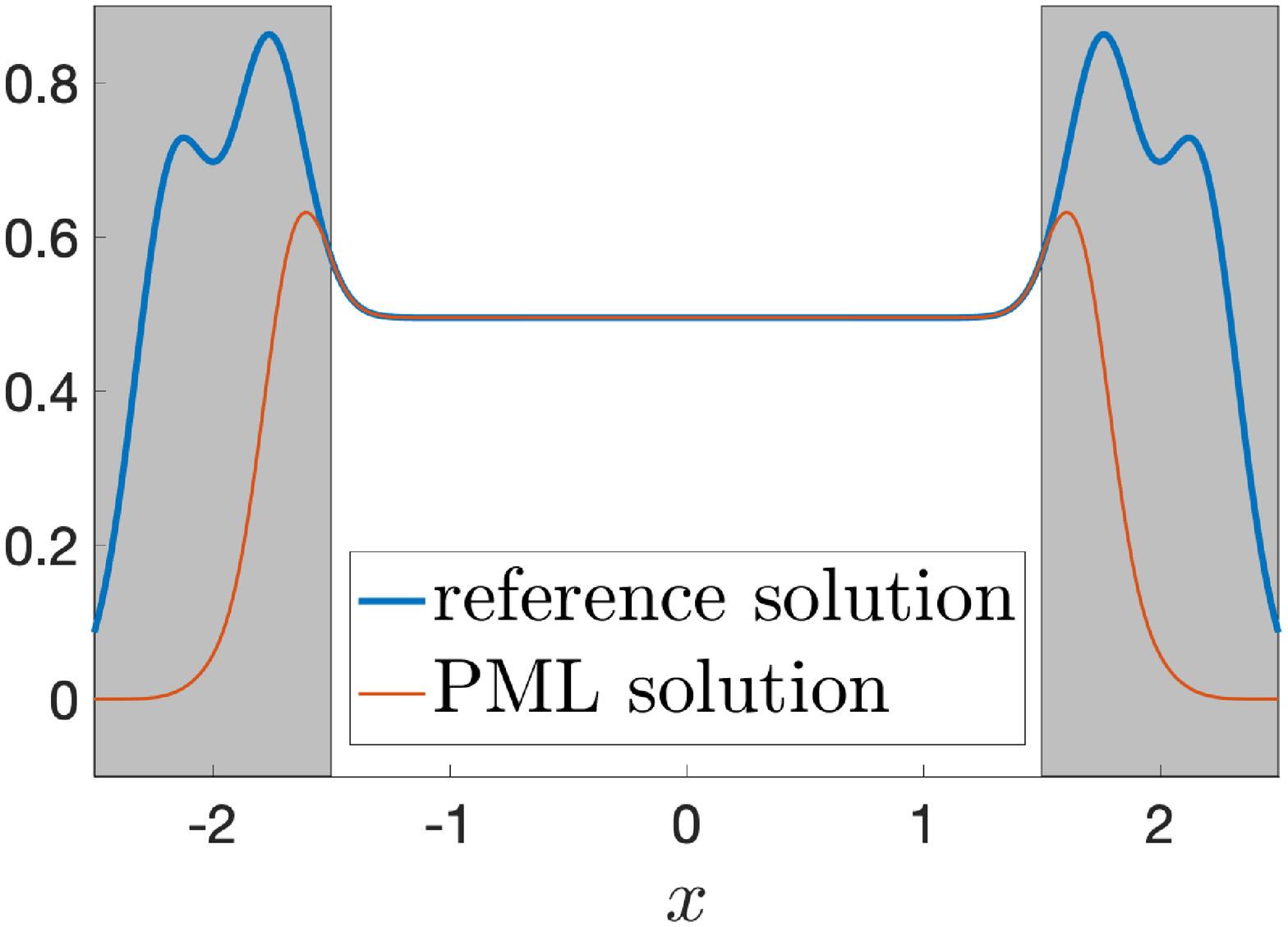}
 \caption{$\delta=0.1$}
 \label{fig:ex5p1solutionsT2c}
 \end{subfigure}
	\caption{(Example 5.1) The comparison of the numerical solution and reference solution at time $t=2$. The PML layers are shaded in light grey.}
	\label{fig_ex5p1solutionsT2}
\end{figure}

 Figure~\ref{fig_ex5p1solutions} shows the evolution of the numerical and reference solutions by taking $T=3, h=2^{-8}, \tau=\frac{1}{12000}$ and $m=400$ for different $\delta=0.5,0.2,0.1$. The reference solutions are obtained by the pseudo-spectral method on a sufficiently large truncated domain. Figure~\ref{fig_ex5p1solutionsT2} shows the PML solution and reference solution at time $t=2$, which indicates the waves decay exponentially in the PML media.

\begin{figure}[htbp]
\centering
\begin{subfigure}{.32\textwidth}
 \centering
	\includegraphics[width=1\textwidth]{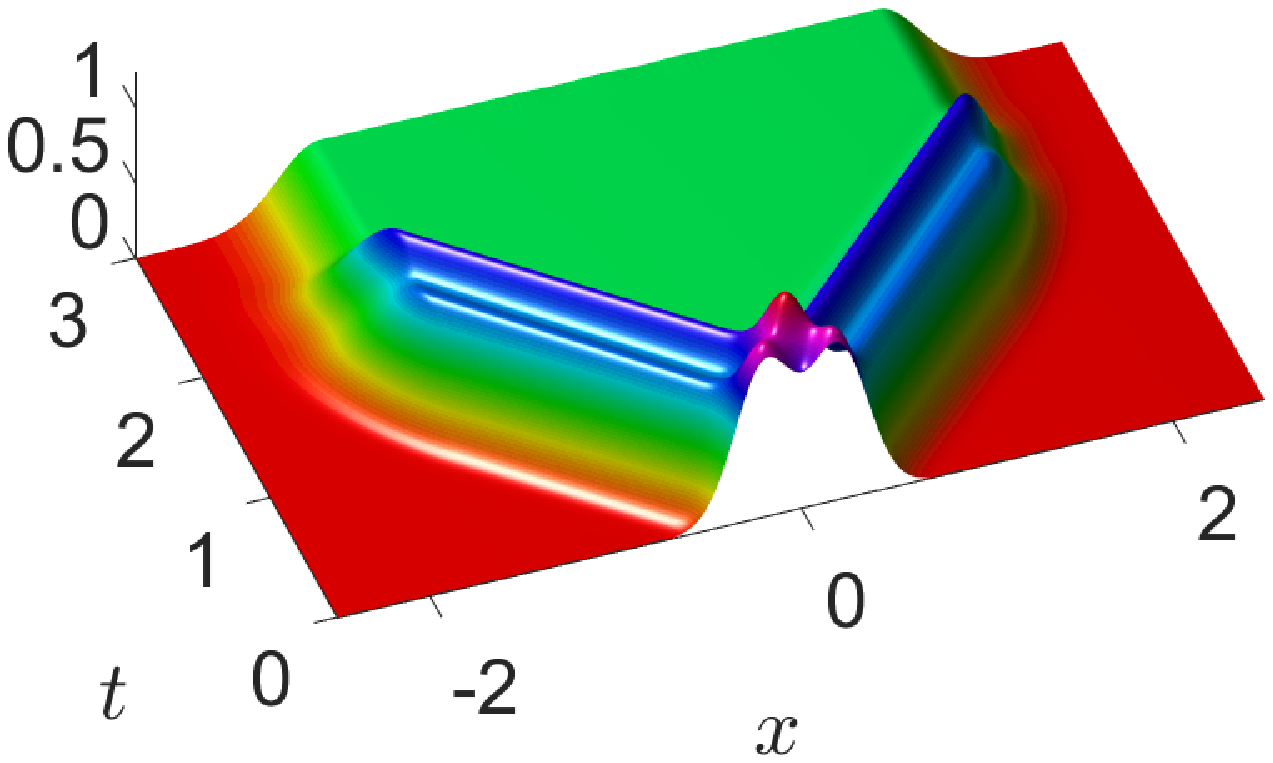}
 \caption{numerical solution $\delta=0.5$}
 \label{fig:ex5p1solutionsa}
\end{subfigure}
\begin{subfigure}{.32\textwidth}
 \centering
	\includegraphics[width=1\textwidth]{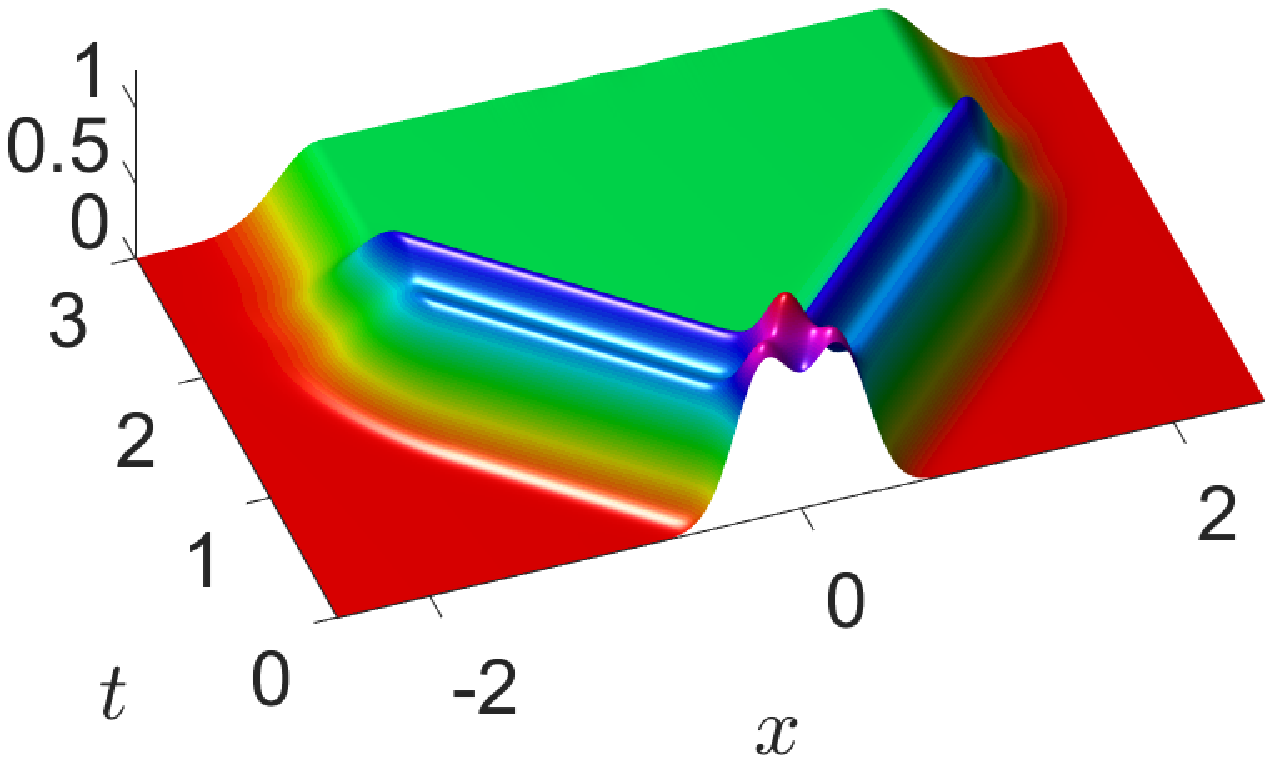}
 \caption{numerical solution $\delta=0.2$}
 \label{fig:ex5p1solutionsb}
\end{subfigure}
\begin{subfigure}{.32\textwidth}
 \centering
	\includegraphics[width=1\textwidth]{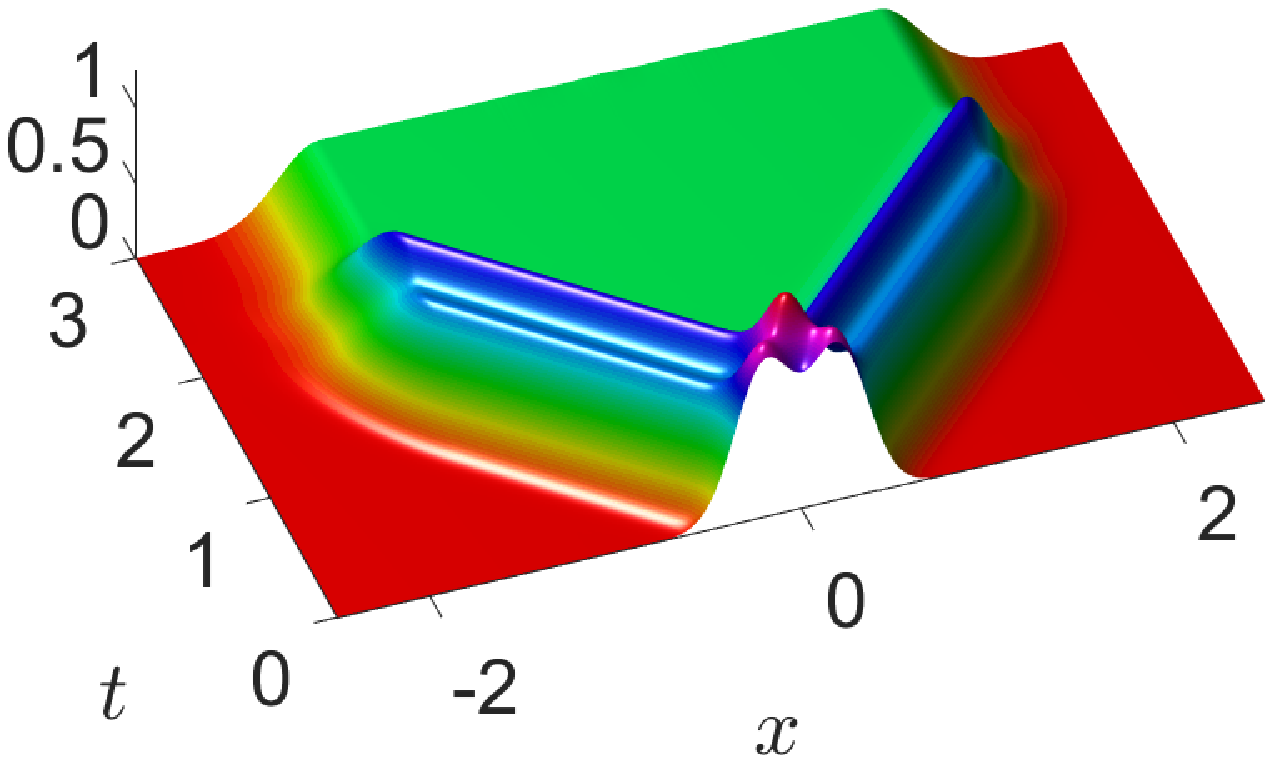}
 \caption{numerical solution $\delta=0.1$}
 \label{fig:ex5p1solutionsc}
 \end{subfigure}
 
 \begin{subfigure}{.32\textwidth}
 \centering
	\includegraphics[width=1\textwidth]{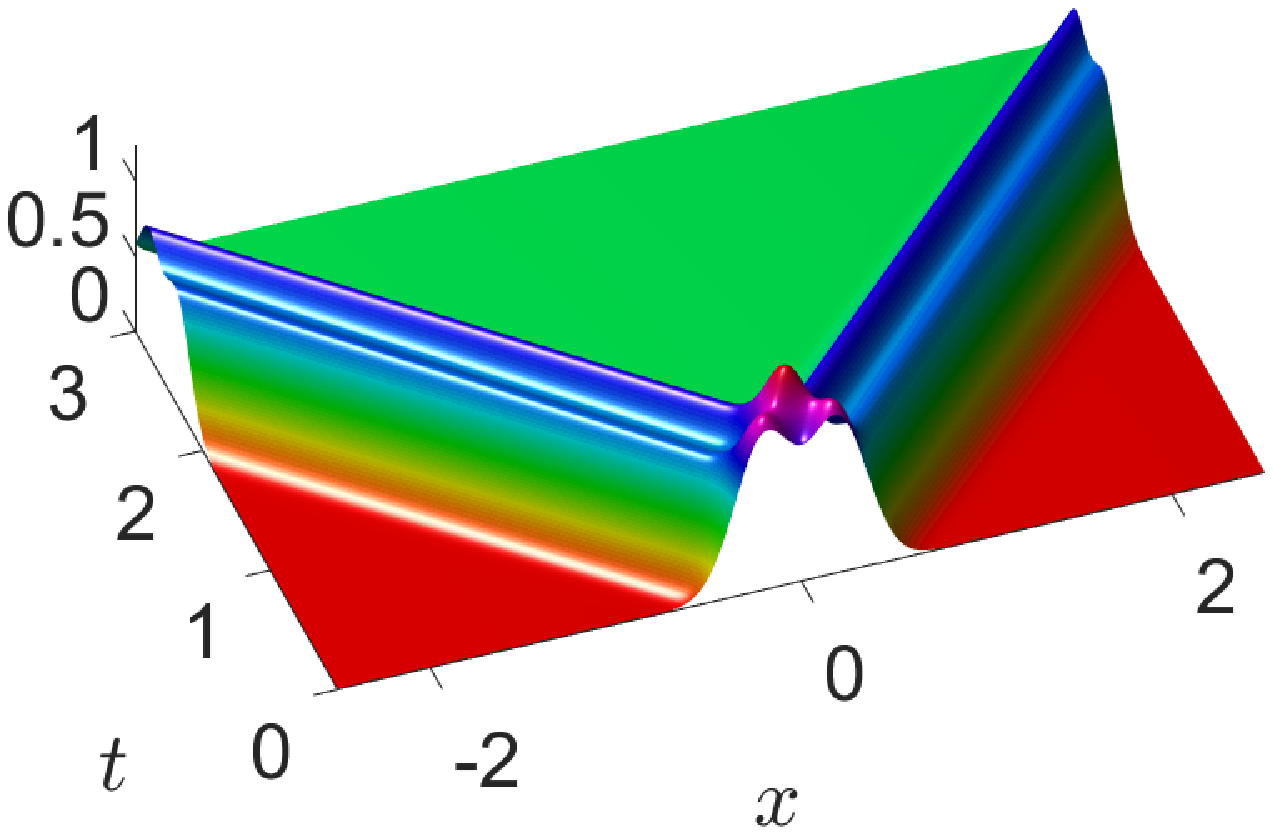}
 \caption{reference solution $\delta=0.5$}
 \label{fig:ex5p1solutionsd}
\end{subfigure}
\begin{subfigure}{.32\textwidth}
 \centering
	\includegraphics[width=1\textwidth]{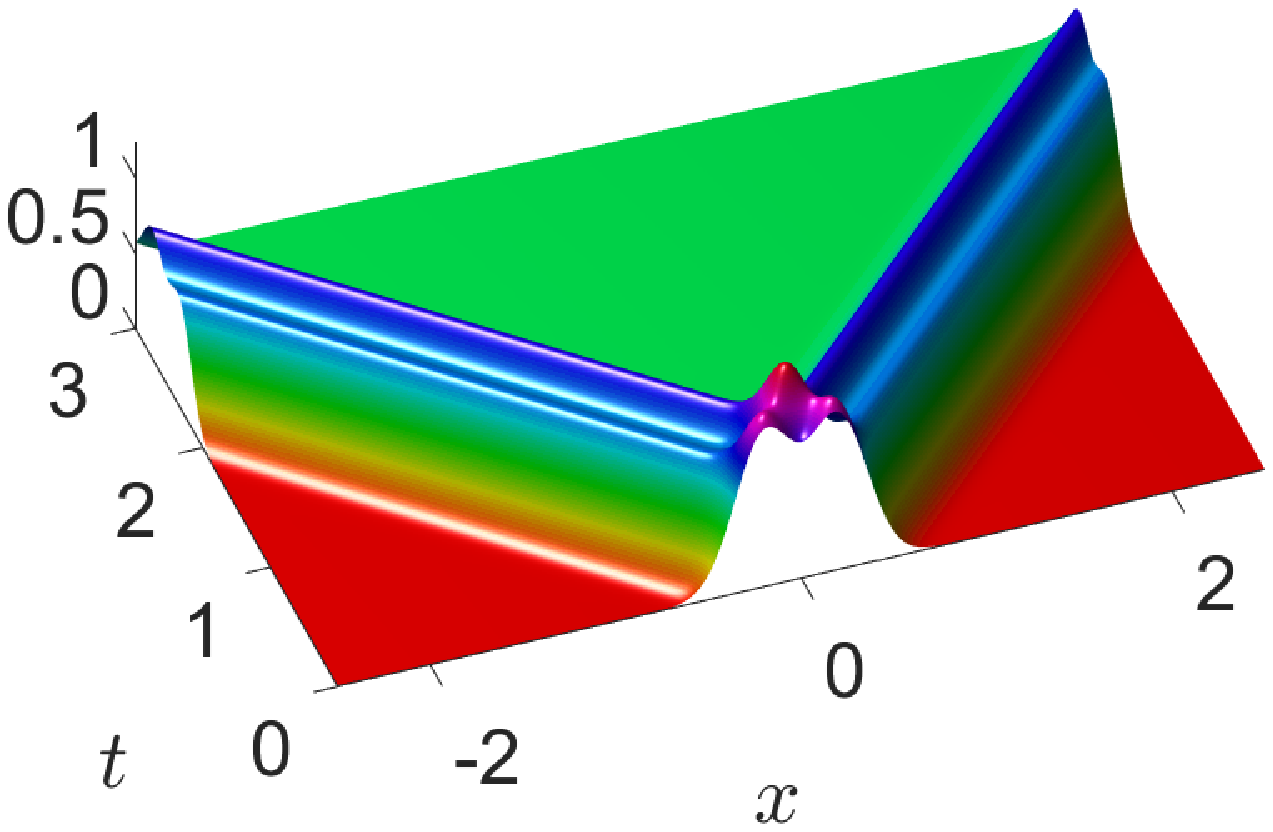}
 \caption{reference solution $\delta=0.2$}
 \label{fig:ex5p1solutionse}
\end{subfigure}
\begin{subfigure}{.32\textwidth}
 \centering
	\includegraphics[width=1\textwidth]{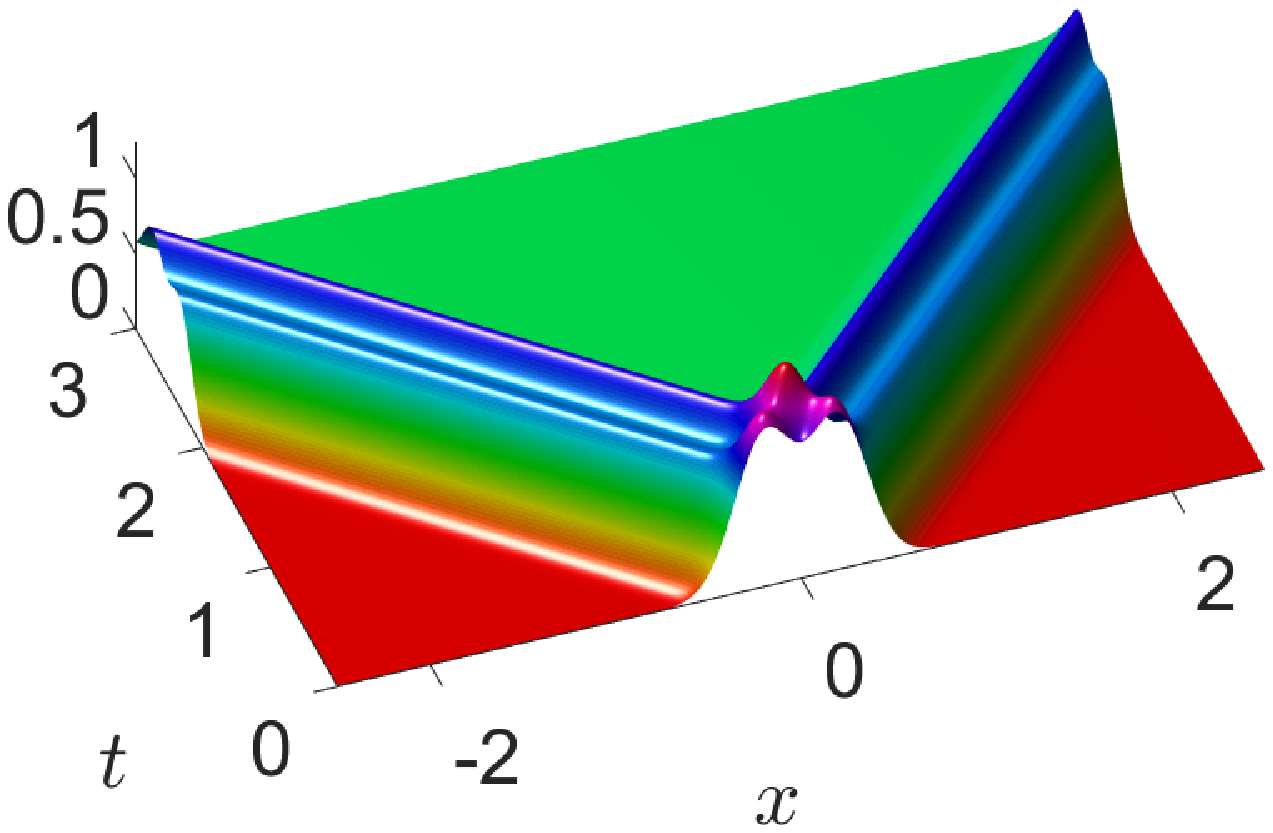}
 \caption{reference solution $\delta=0.1$}
 \label{fig:ex5p1solutionsf}
 \end{subfigure}
	\caption{(Example 5.1) Numerical solutions for $\delta=0.5,0.2,0.1$ and reference solutions up to $T=3$. The numerical solutions are obtained by taking $h=2^{-8}$ and $m=400$..}
	\label{fig_ex5p1solutions}
\end{figure}

\begin{table}
\centering
\begin{tabular}{c|cc|cc|cc}
\hline
$h$	&	$\delta=0.5$  & order & $\delta=0.2$ & order &  $\delta=0.1$ & order\\ \hline
$2^{-4}$	&  8.39e-03   & -- & 7.62e-03 & -- & 7.91e-03 & --\\
$2^{-5}$	& 1.36e-03 &  2.62 & 1.19e-03 & 2.67 & 1.31e-03 & 2.59\\
$2^{-6}$ &  3.16e-04 & 2.11 & 2.68e-04 & 2.15 & 2.87e-04 & 2.19\\
$2^{-7}$ & 7.79e-05 & 2.02 & 6.32e-05 & 2.09 & 6.73e-05 & 2.09\\\hline
\end{tabular}
\caption{(Example 5.1) Errors $e_h$ and convergence orders at $t=2$.}\label{tab_5p1_l2}
\end{table}

Table~\ref{tab_5p1_l2} shows the errors $e_h$ at time $t=2$ and the spatial convergence rates for various horizons $\delta=0.5,0.2,0.1$.
We now investigate if numerical solutions ${\tilde q}$ of \eqref{eq_truPML1}--\eqref{eq_truPML3} converge to the correct solution $u$ of the corresponding local problem \eqref{eq_localPML1}--\eqref{eq_localPML2} as the horizon $\delta$ goes to zero. The local wave equations with PML modifications in \cite{basu2008perfectly} are given by
\begin{align}
\frac{\partial^2 u}{\partial t^2} + z\sigma(x)\frac{\partial u}{\partial t} - \frac{\partial v}{\partial x}&= 0,\label{eq_localPML1}\\
\frac{\partial v}{\partial t} - \frac{\partial^2 u}{\partial t\partial x} + z\sigma(x)v & = 0.\label{eq_localPML2}
\end{align}
In the simulation, we fix the ratio $\delta/h=M$ with $M=1,2,3$. Table~\ref{tab_5p1_dh_l2} shows the errors $e_\delta$ at $t=2$ and the second-order convergence rates, which is consistent to the analysis in \cite{TianDu}.

\begin{table}
\centering
\begin{tabular}{c|cc|cc|cc}
\hline
$h$	&	$\delta=h$  & order & $\delta=2h$ & order &  $\delta=3h$ & order\\ \hline
$2^{-4}$	&  2.63e-02   &-- & 2.62e-02 & -- & 2.62e-03 & --\\
$2^{-5}$	&  5.72e-03 &  2.20 & 5.71e-03 & 2.20 & 5.70e-03 & 2.20\\
$2^{-6}$ &  1.35e-03  & 2.08 & 1.34e-03 & 2.09 & 1.35e-04 & 2.09\\
$2^{-7}$ & 3.15e-04 & 2.09 & 3.14e-04 & 2.10 & 3.13e-05 & 2.09\\\hline
\end{tabular}
\caption{(Example 5.1) Errors $e_\delta$ and $\delta$-convergence orders between the numerical solutions and exact solutions of local problem \eqref{eq_localPML1}--\eqref{eq_localPML2} by vanishing $\delta$ and $h$ simultaneously at $t=2$.}\label{tab_5p1_dh_l2}
\end{table}

\bigskip 

\noindent\textbf{Example 5.2.} In this example, we take the source $f(x,t)\equiv0$ and the initial values
$$\psi_0(x) = e^{-25(x-0.2)^2} + e^{-25(x+0.2)^2}, \quad \psi_1(x)  = 50xe^{-25x^2},$$ 
and consider the Gaussian kernel in the form of 
\begin{align}
\gamma(y-x,\frac{x+y}{2}) = \frac{4}{\delta^3}\sqrt{\frac{10^3}{\pi}} e^{-10\frac{(x-y)^2}{\delta^2}}. \label{eq_gaukernel}
\end{align}
In the simulations, we set $l=2,d_p=2$, the PML coefficient $z=10$, and take Talbot's contour parameters $\mu = |z|$ and $\nu = 1$. The reason for choosing these parameters is discussed in Appendix~\ref{sec_Ap}.

\begin{figure}[htbp]
\centering
\begin{subfigure}{.32\textwidth}
 \centering
	\includegraphics[width=1\textwidth]{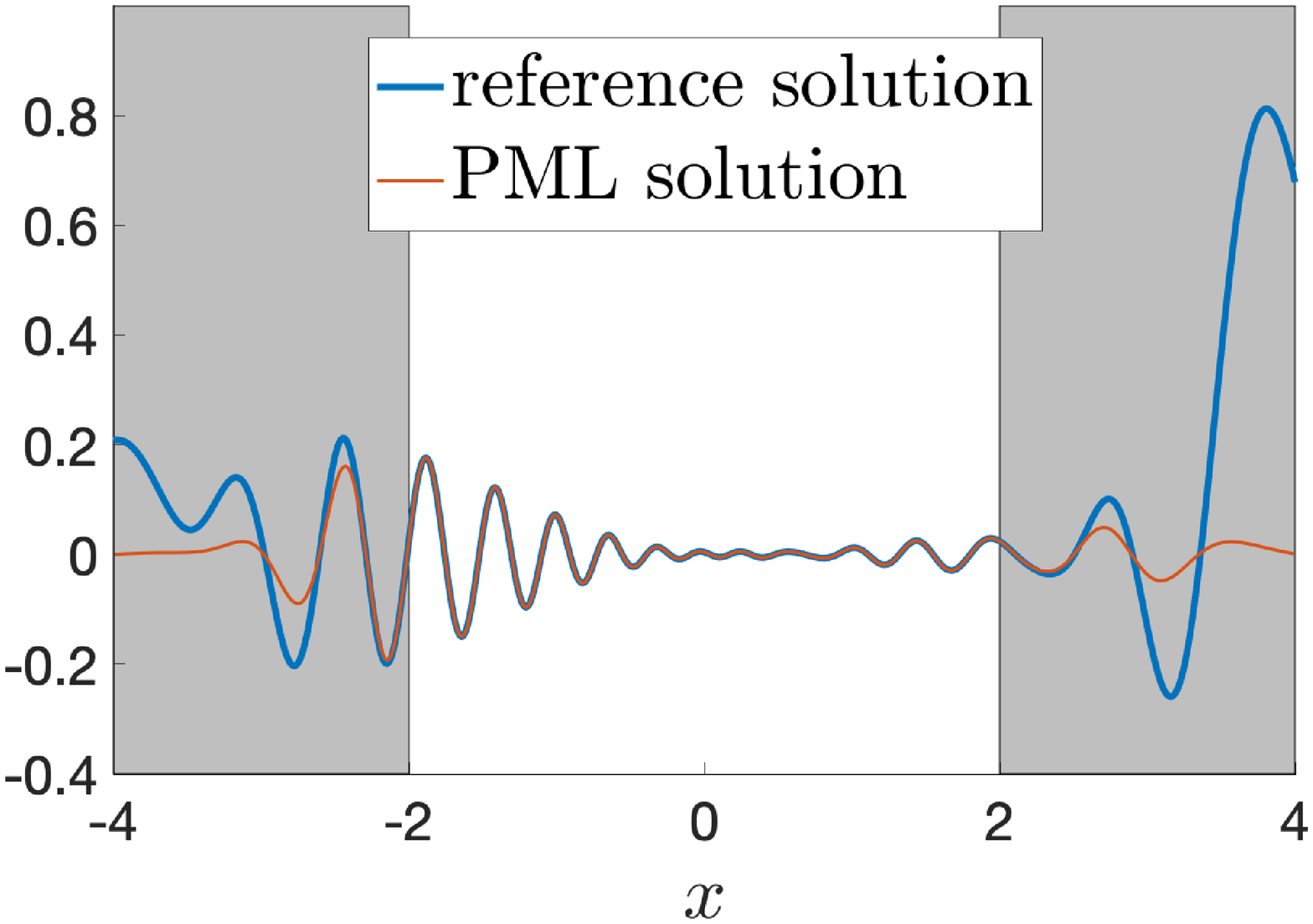}
 \caption{$\delta=0.5$}
 \label{fig:ex5p2solutionsT2a}
\end{subfigure}
\begin{subfigure}{.32\textwidth}
 \centering
	\includegraphics[width=1\textwidth]{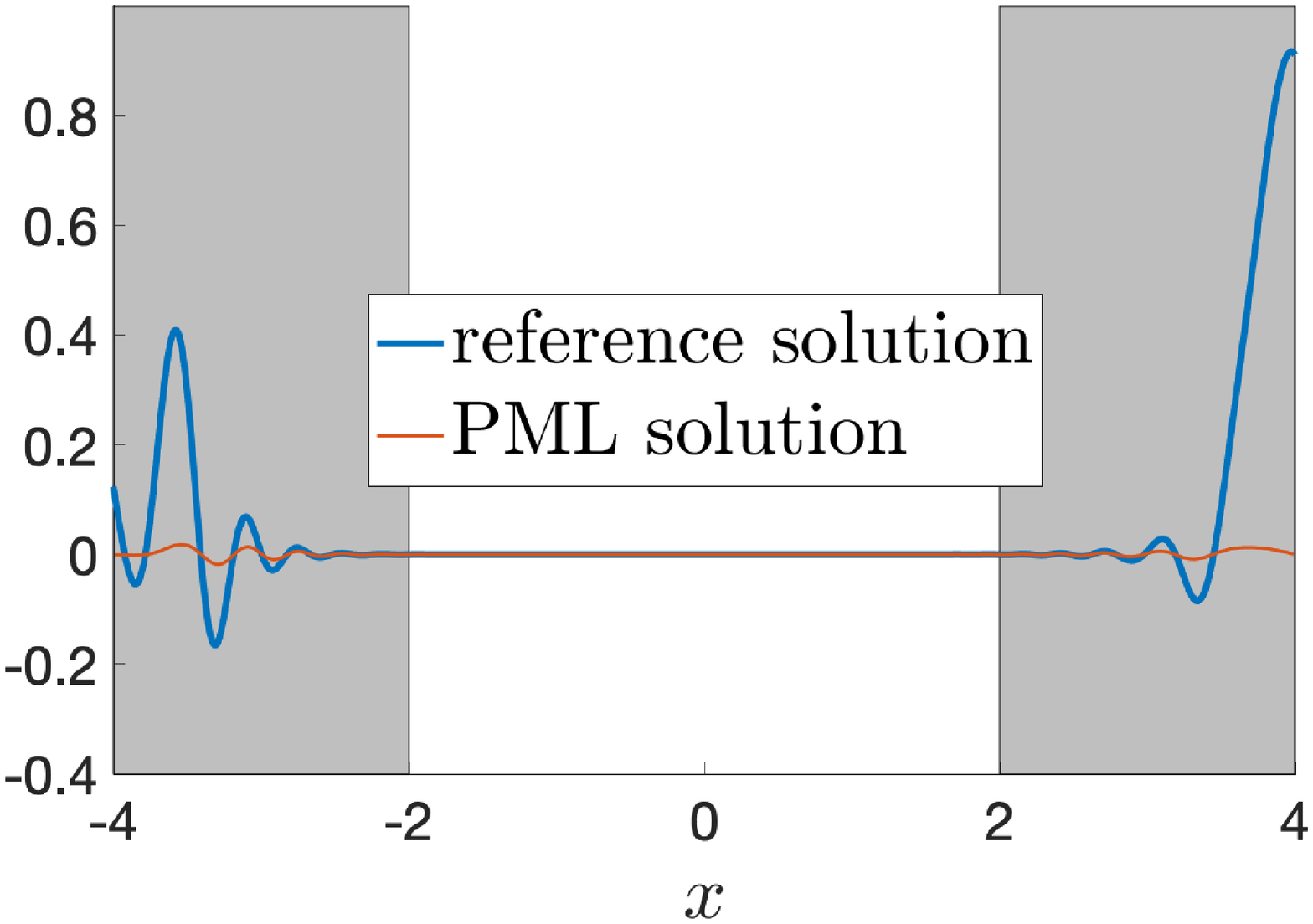}
 \caption{$\delta=0.2$}
 \label{fig:ex5p2solutionsT2b}
\end{subfigure}
\begin{subfigure}{.32\textwidth}
 \centering
	\includegraphics[width=1\textwidth]{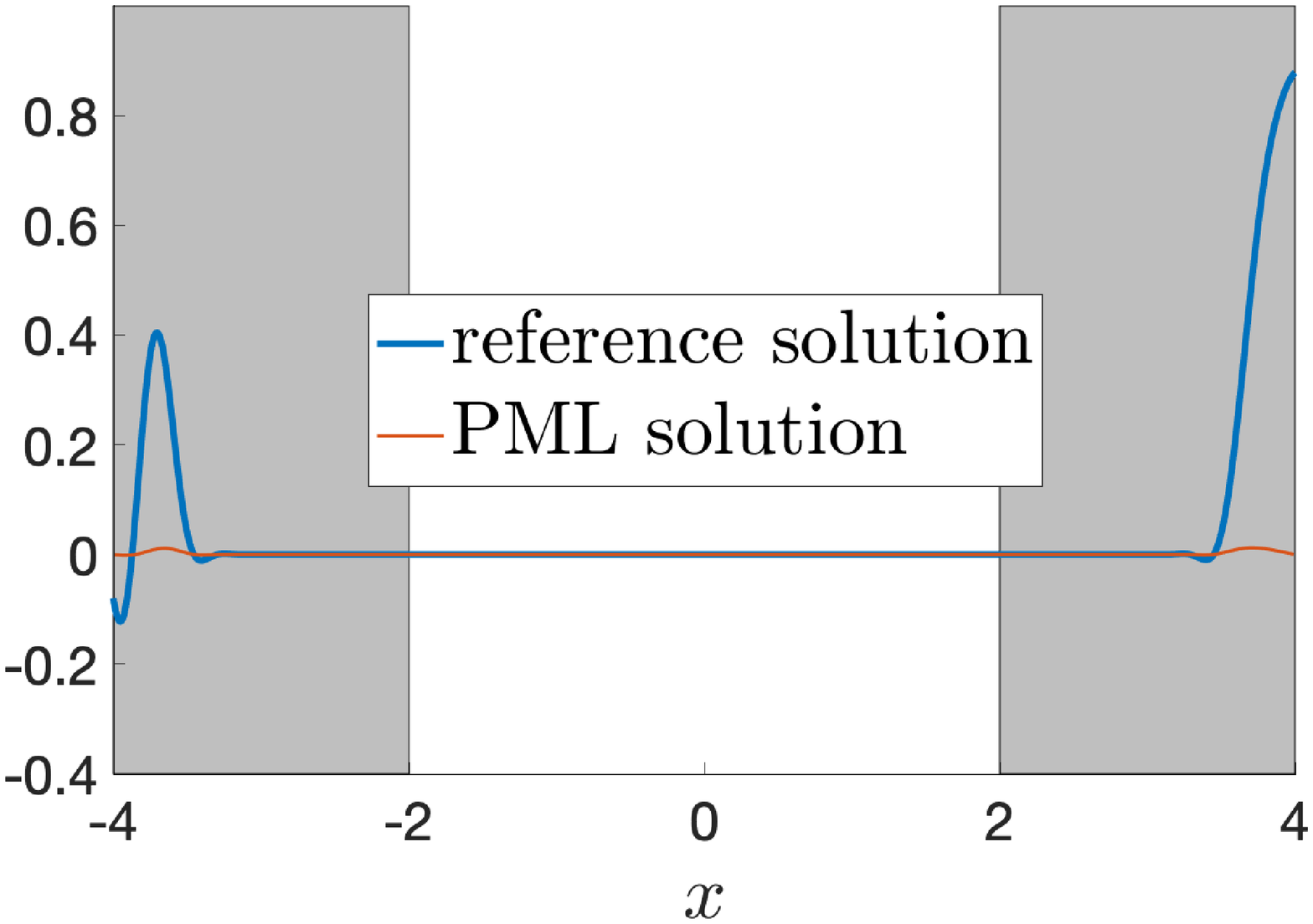}
 \caption{$\delta=0.1$}
 \label{fig:ex5p2solutionsT2c}
 \end{subfigure}
	\caption{(Example 5.2) The comparison of the numerical solution and reference solution at time $t=4$. The PML layers for numerical solutions are shaded in light grey.}
	\label{fig_ex5p2solutionsT2}
\end{figure}

Figure~\ref{fig_ex5p2solutions} plots the evolution of numerical and the reference solutions by taking $T=4$, $h=2^{-8}$, $\tau=\frac{1}{12000}$ and $m=800$ for different $\delta=0.5,0.2,0.1$. Figure~\ref{fig_ex5p2solutionsT2} plots the PML solution and the reference solution at time $t=2$, which indicates the waves decay exponentially in the PML media.

\begin{figure}[htbp]
\centering
\begin{subfigure}{.32\textwidth}
 \centering
	\includegraphics[width=1\textwidth]{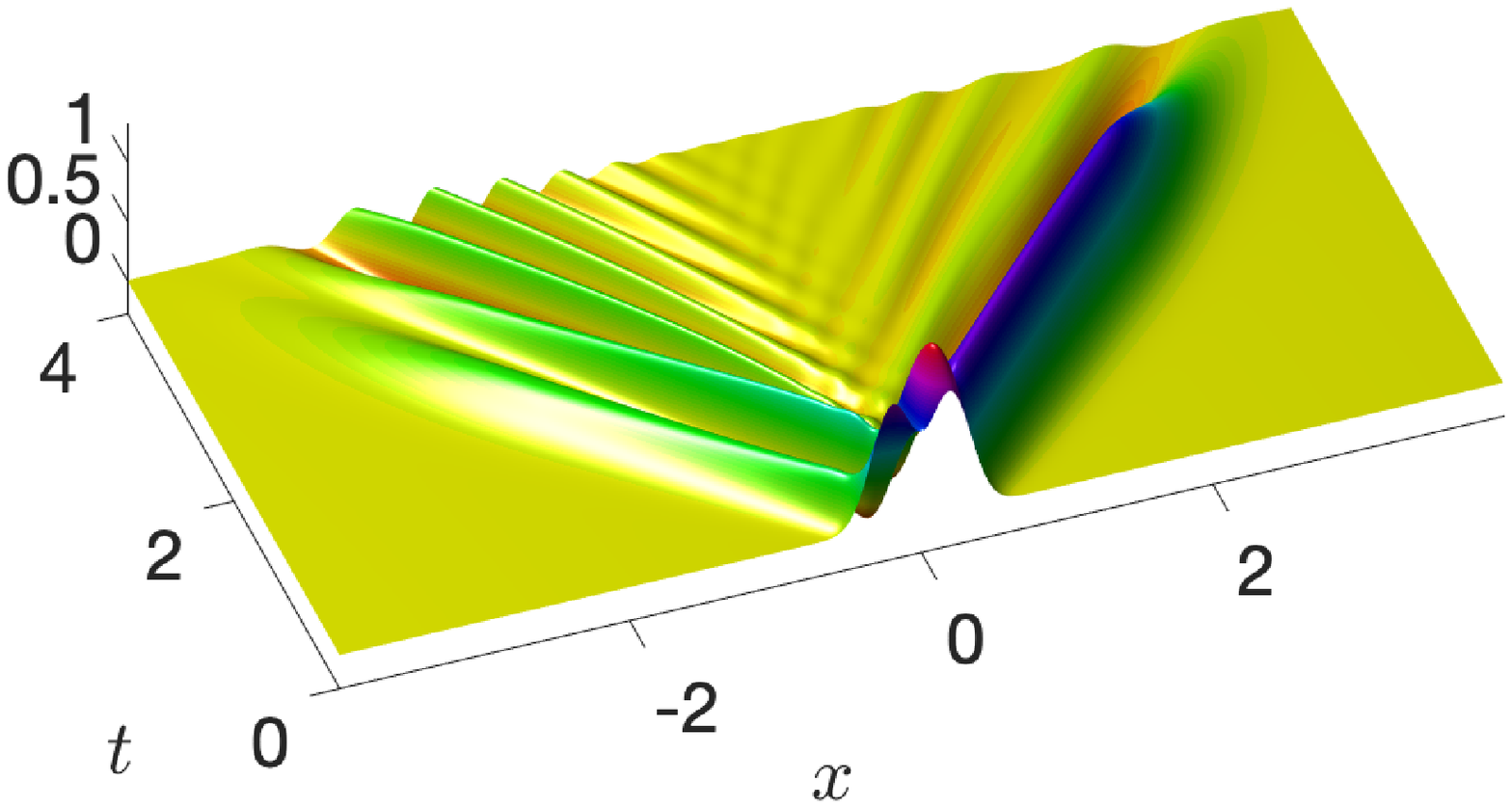}
 \caption{numerical solution $\delta=0.5$}
 \label{fig:ex5p2solutionsa}
\end{subfigure}
\begin{subfigure}{.32\textwidth}
 \centering
	\includegraphics[width=1\textwidth]{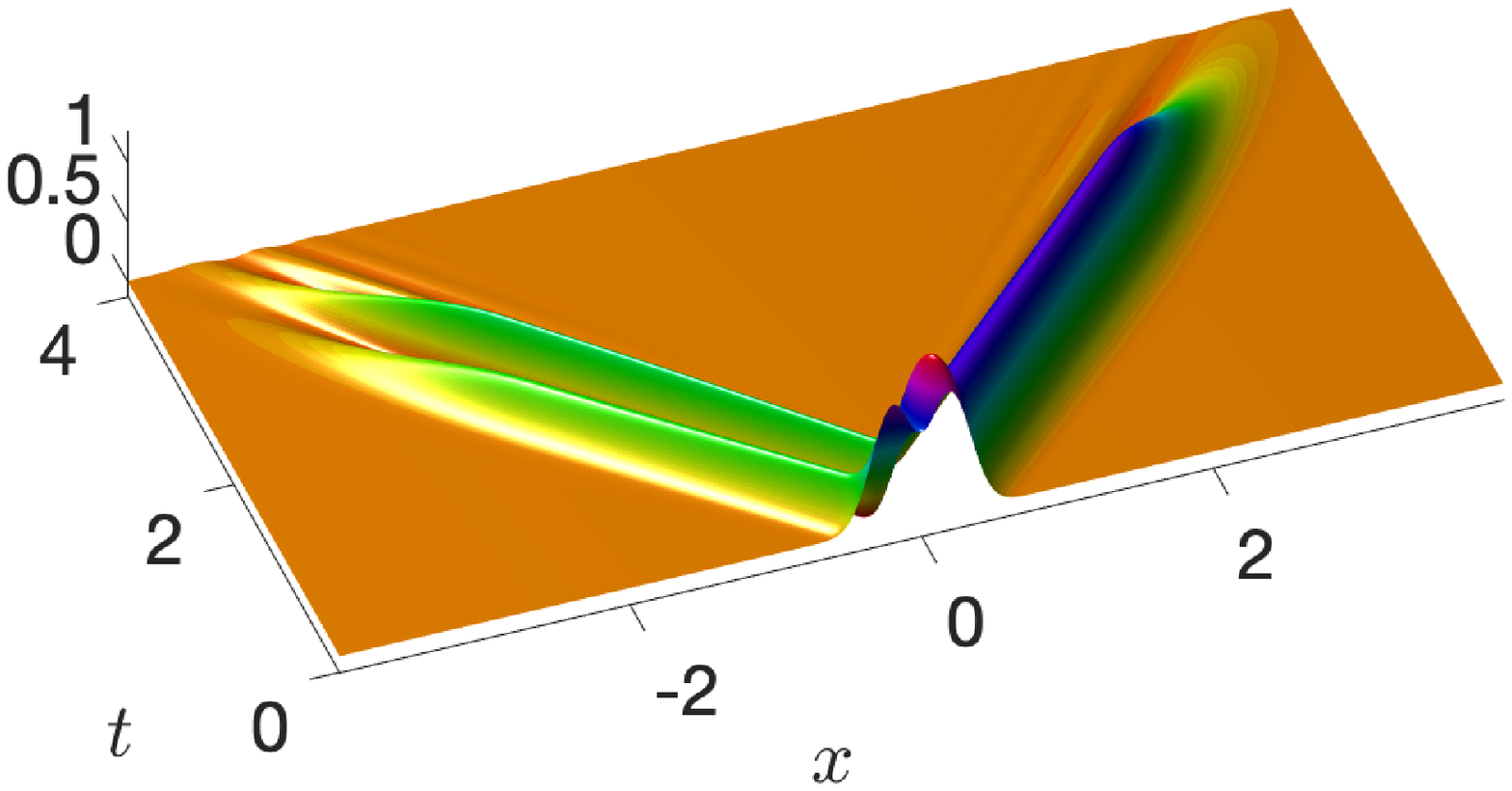}
 \caption{numerical solution $\delta=0.2$}
 \label{fig:ex5p2solutionsb}
\end{subfigure}
\begin{subfigure}{.32\textwidth}
 \centering
	\includegraphics[width=1\textwidth]{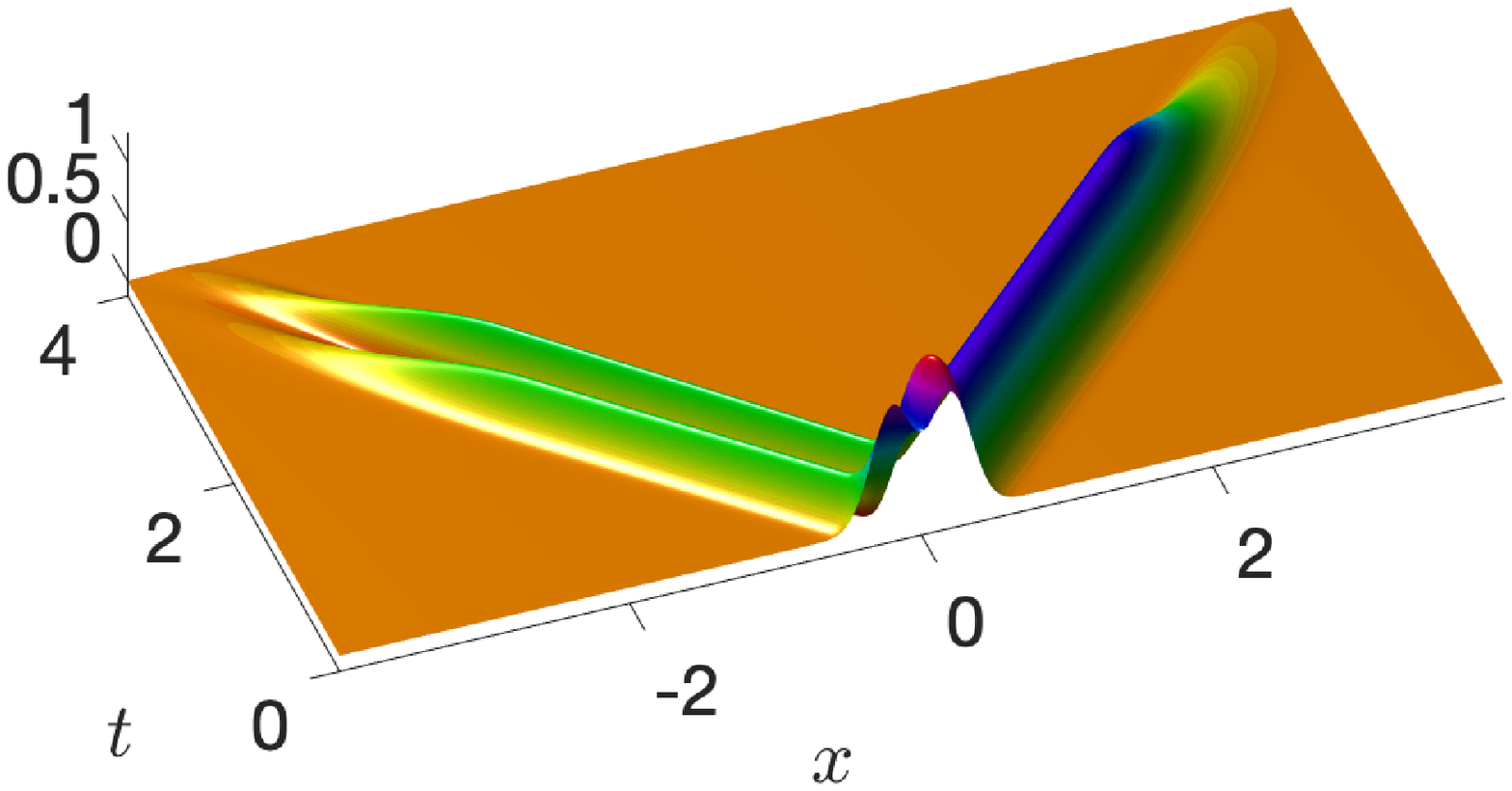}
 \caption{numerical solution $\delta=0.1$}
 \label{fig:ex5p2solutionsc}
 \end{subfigure}
 
 \begin{subfigure}{.32\textwidth}
 \centering
	\includegraphics[width=1\textwidth]{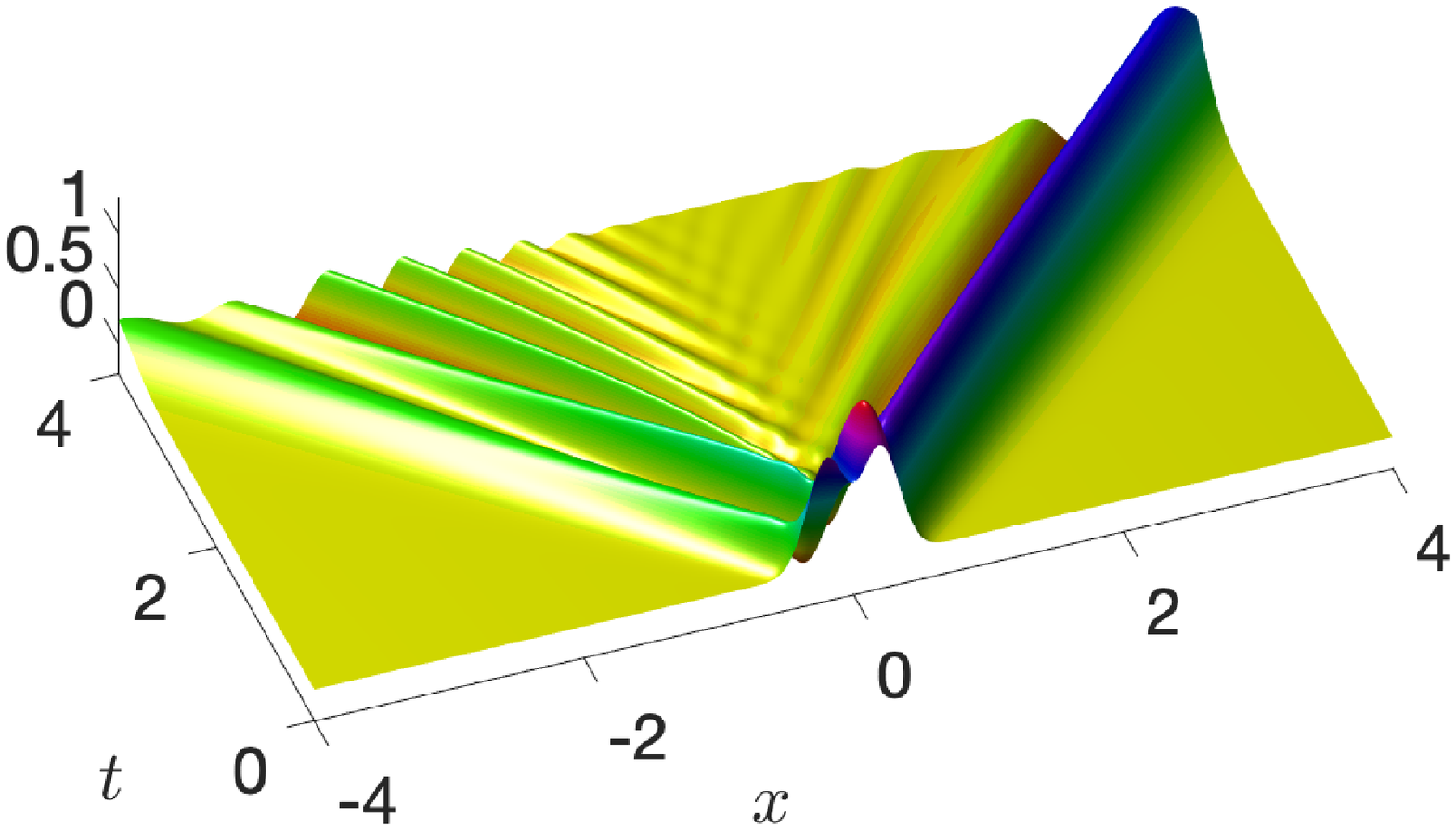}
 \caption{reference solution $\delta=0.5$}
 \label{fig:ex5p2solutionsd}
\end{subfigure}
\begin{subfigure}{.32\textwidth}
 \centering
	\includegraphics[width=1\textwidth]{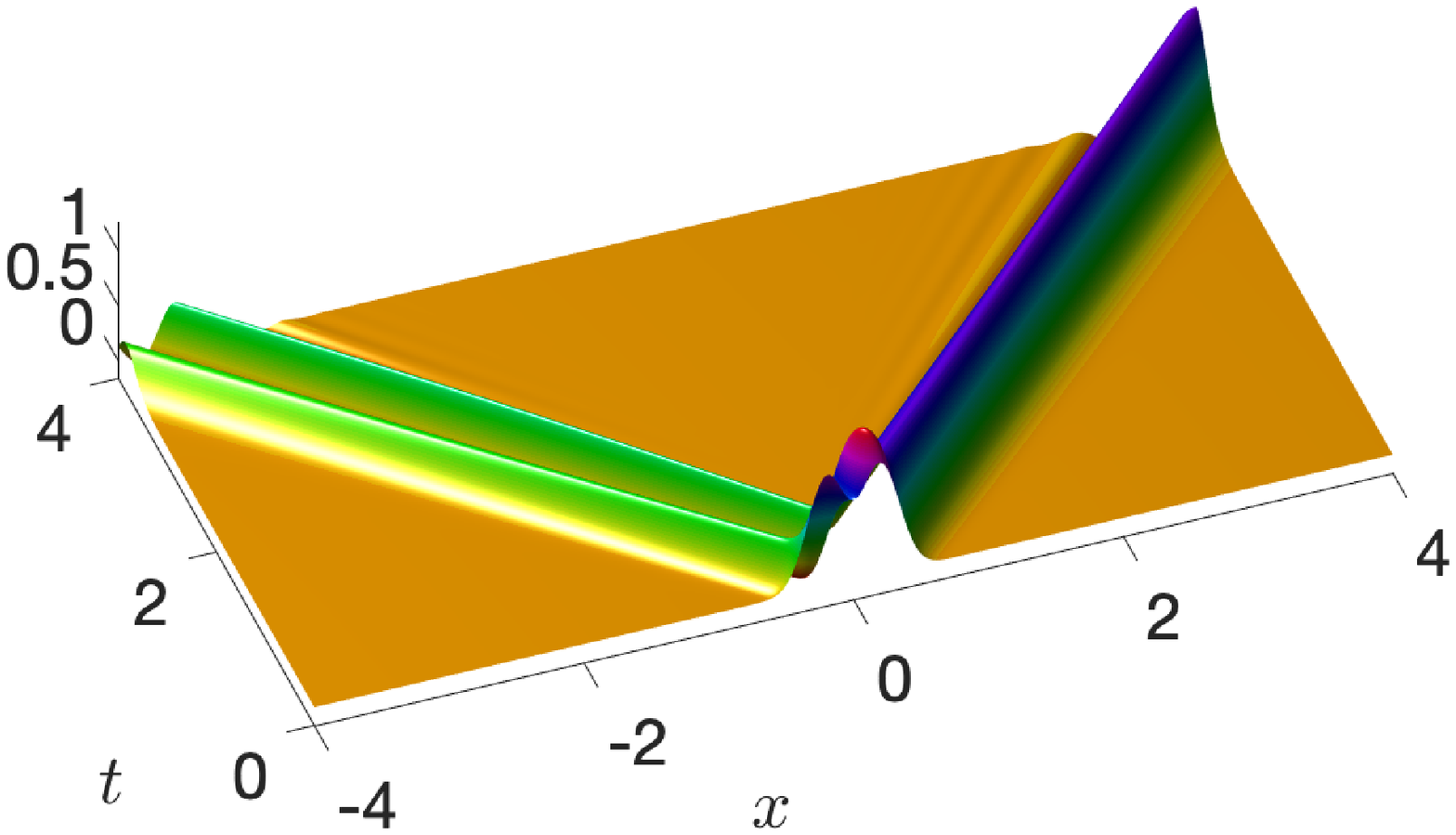}
 \caption{reference solution $\delta=0.2$}
 \label{fig:ex5p2solutionse}
\end{subfigure}
\begin{subfigure}{.32\textwidth}
 \centering
	\includegraphics[width=1\textwidth]{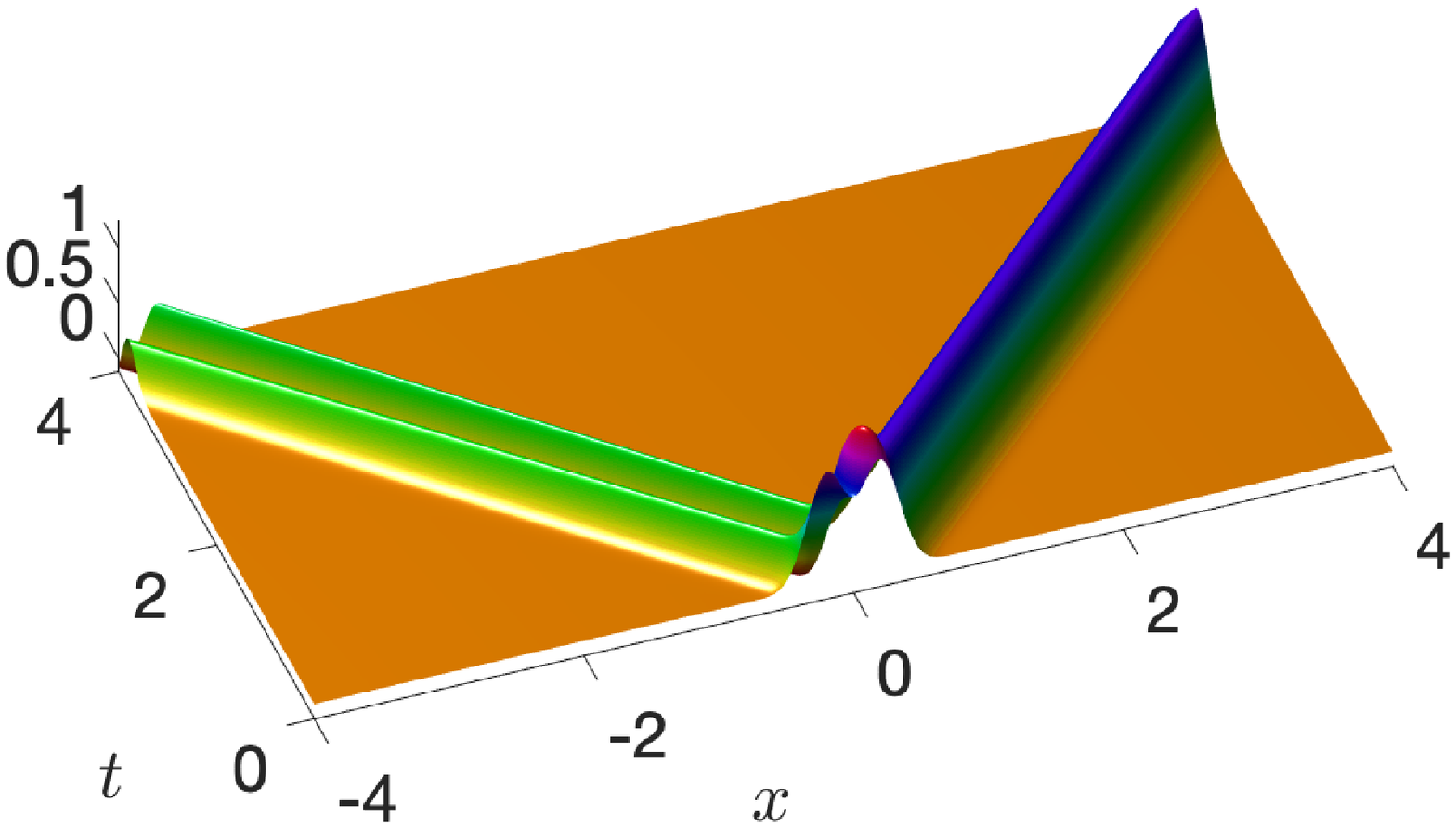}
 \caption{reference solution $\delta=0.1$}
 \label{fig:ex5p2solutionsf}
 \end{subfigure}
	\caption{(Example 5.2) Numerical solutions for $\delta=0.5,0.2,0.1$ and reference solutions up to $T=4$. The numerical solutions are obtained by taking $h=2^{-8}$ and $m=400$.}
	\label{fig_ex5p2solutions}
\end{figure}

Table~\ref{tab_5p2_l2} shows the errors $e_h$ and the spatial convergence order at $t=4$ by refining $h$ for $\tau=\frac{1}{12000}$ and $m=400$. Table~\ref{tab_5p2_dh_l2} shows the errors $e_\delta$ and $\delta$-convergence rates of second-order with the corresponding local PML equations~\eqref{eq_localPML1}--\eqref{eq_localPML2}.
\begin{table}
\centering
\begin{tabular}{c|cc|cc|cc}
\hline
$h$	&	$\delta=0.5$  & order & $\delta=0.2$ & order &  $\delta=0.1$ & order\\ \hline
$2^{-4}$	&  4.07e-02   & -- & 3.18e-02 & -- & 1.50e-02 & --\\
$2^{-5}$	& 1.04e-02 &  1.96 & 7.19e-03 & 2.15 & 2.38e-03 & 2.65\\
$2^{-6}$ &  2.63e-03 & 1.99 & 1.74e-03 & 2.05 & 5.62e-04 & 2.08\\
$2^{-7}$ & 6.57e-04 & 2.00 & 4.30e-04 & 2.02 & 1.38e-05 & 2.02\\\hline
\end{tabular}
\caption{(Example 5.2) Errors $e_h$ and convergence orders at $t=4$.}\label{tab_5p2_l2}
\end{table}

\begin{table}
\centering
\begin{tabular}{c|cc|cc|cc}
\hline
$h$	&	$\delta=h$  & order & $\delta=2h$ & order &  $\delta=3h$ & order\\ \hline
$2^{-4}$	& 4.47e-03   & -- & 5.73e-03  & -- & 1.30e-02  & --\\
$2^{-5}$	&  8.29e-04 &  2.43 & 9.89e-04  & 2.53 & 1.82e-03 & 2.83\\
$2^{-6}$ &  1.92e-04  & 2.11 & 2.24e-04 & 2.14 & 3.91e-04  & 2.22\\
$2^{-7}$ & 4.43e-05 & 2.11 & 5.17e-05 & 2.11 & 9.16e-05& 2.09\\\hline
\end{tabular}
\caption{(Example 5.2) Errors $e_\delta$ and $\delta$-convergence orders between the numerical solutions and exact solutions of local problem \eqref{eq_localPML1}--\eqref{eq_localPML2} by vanishing $\delta$ and $h$ simultaneously at $t=4$.}\label{tab_5p2_dh_l2}
\end{table}

\bigskip 
\noindent\textbf{Example 5.3.} Here we use the same source term and the initial values as Example 5.2, and consider the following spatially inhomogeneous kernel
\begin{align}
\gamma(y-x,\frac{x+y}{2}) = \frac{\omega(\frac{x+y}{2})}{\zeta^3(\frac{x+y}{2})} H\Big(\frac{y-x}{\zeta(\frac{x+y}{2})}\Big),
\end{align}
where 
\begin{align}
\omega(\beta) = 1+e^{-3\beta^2},\ \zeta(\beta) = \delta(2+\tanh(-1.5\beta)),\ H(s) = 4\sqrt{\frac{10^3}{\pi}} e^{-10s^2}.
\end{align}
In the simulations, we take $l=2, d_p=2$, the PML coefficient $z=10$ and $\mu = |z|, \nu = 1$.

\begin{figure}[htbp]
\centering
\begin{subfigure}{.4\textwidth}
 \centering
	\includegraphics[width=1\textwidth]{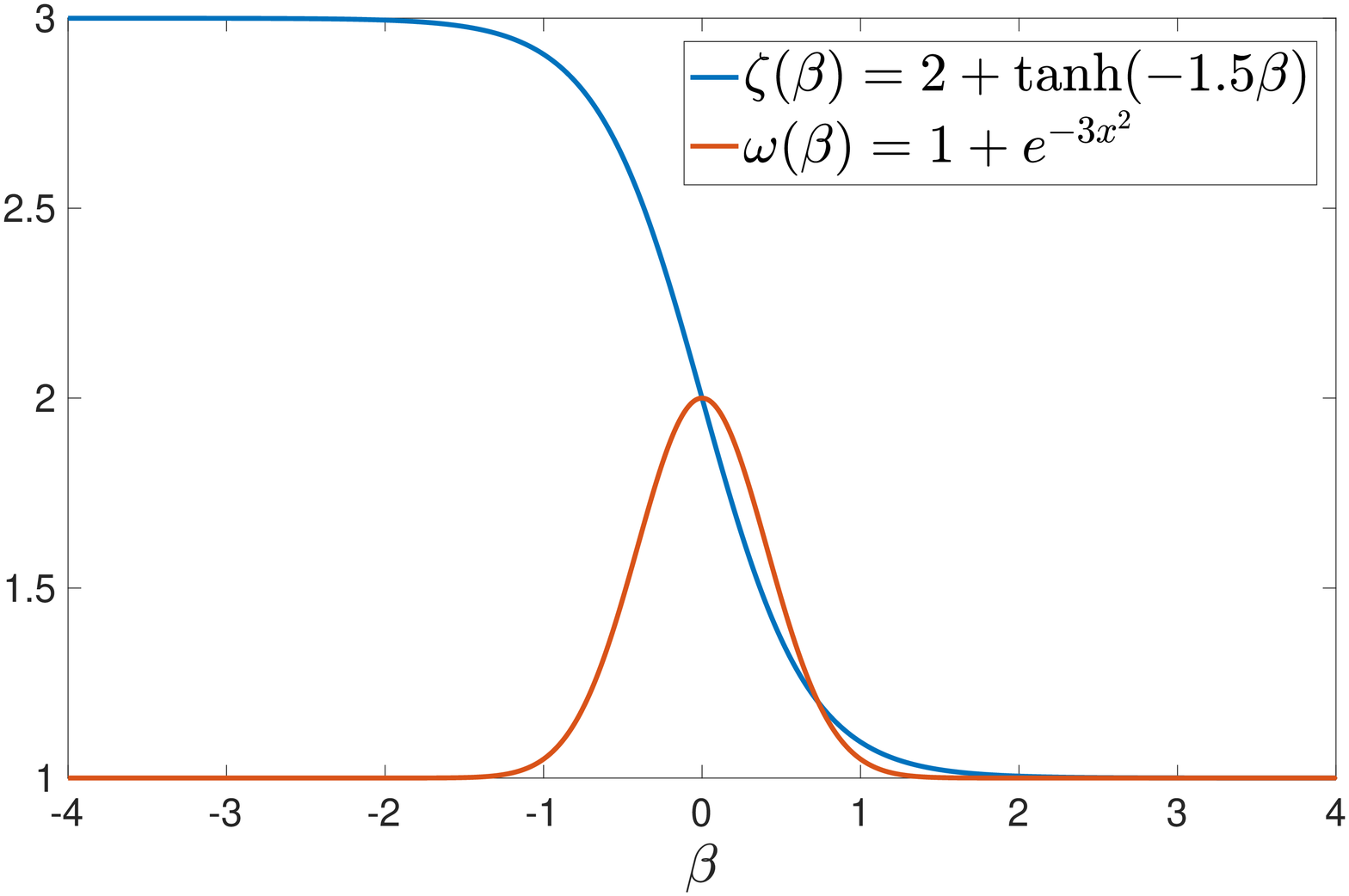}
\end{subfigure}
\begin{subfigure}{.4\textwidth}
 \centering
	\includegraphics[width=1\textwidth]{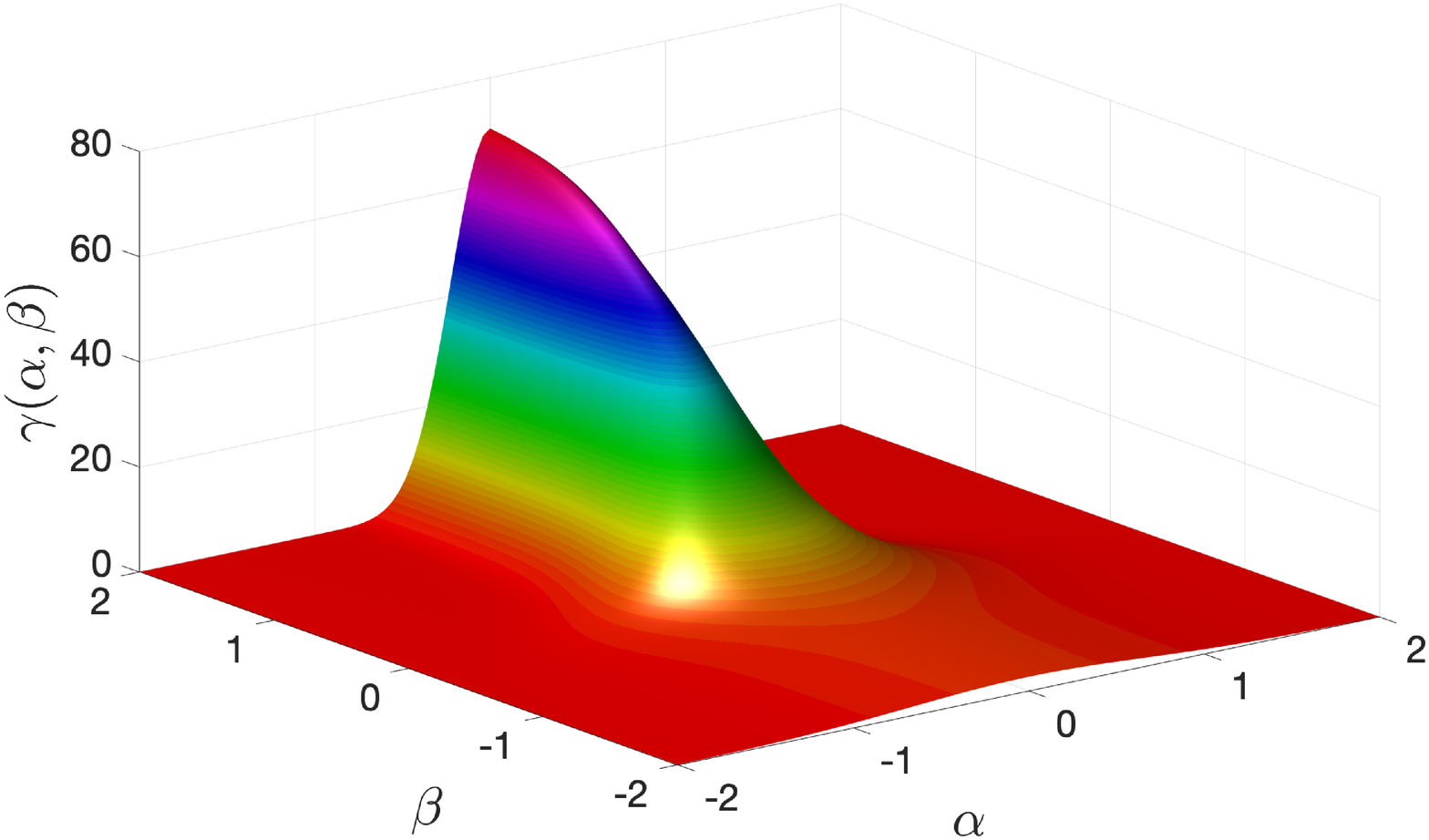}
\end{subfigure}
	\caption{(Example 5.3) Left: the functions $\omega(\beta)$ and $\zeta(\beta)$. Right: the kernel $\gamma(\alpha,\beta)$ for $\delta=1$.}
	\label{fig_ex5p4ker}
\end{figure}

The limiting local wave equation with PML modifications is given by
\begin{align}
\frac{\partial^2 u}{\partial t^2} + z\sigma(x)\frac{\partial u}{\partial t} - \frac{\partial v}{\partial x}&= 0,\label{eq_localinPML1}\\
\frac{\partial v}{\partial t} - \mu(x)\frac{\partial^2 u}{\partial t\partial x} + z\sigma(x)v & = 0.\label{eq_localinPML2}
\end{align}

\begin{figure}[htbp]
\centering
\begin{subfigure}{.32\textwidth}
 \centering
	\includegraphics[width=1\textwidth]{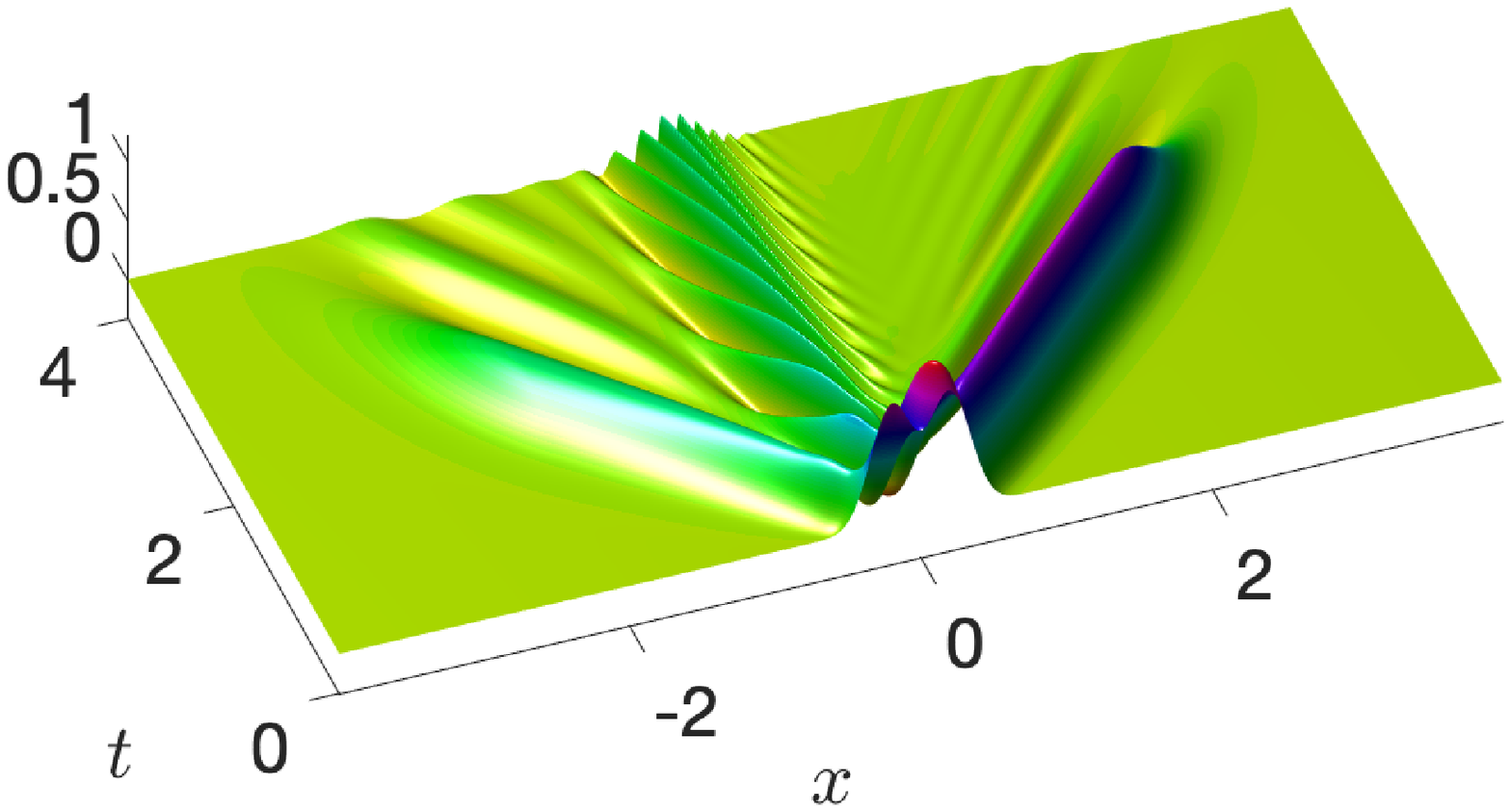}
 \caption{numerical solution $\delta=0.3$}
 \label{fig:ex5p3solutionsa}
\end{subfigure}
\begin{subfigure}{.32\textwidth}
 \centering
	\includegraphics[width=1\textwidth]{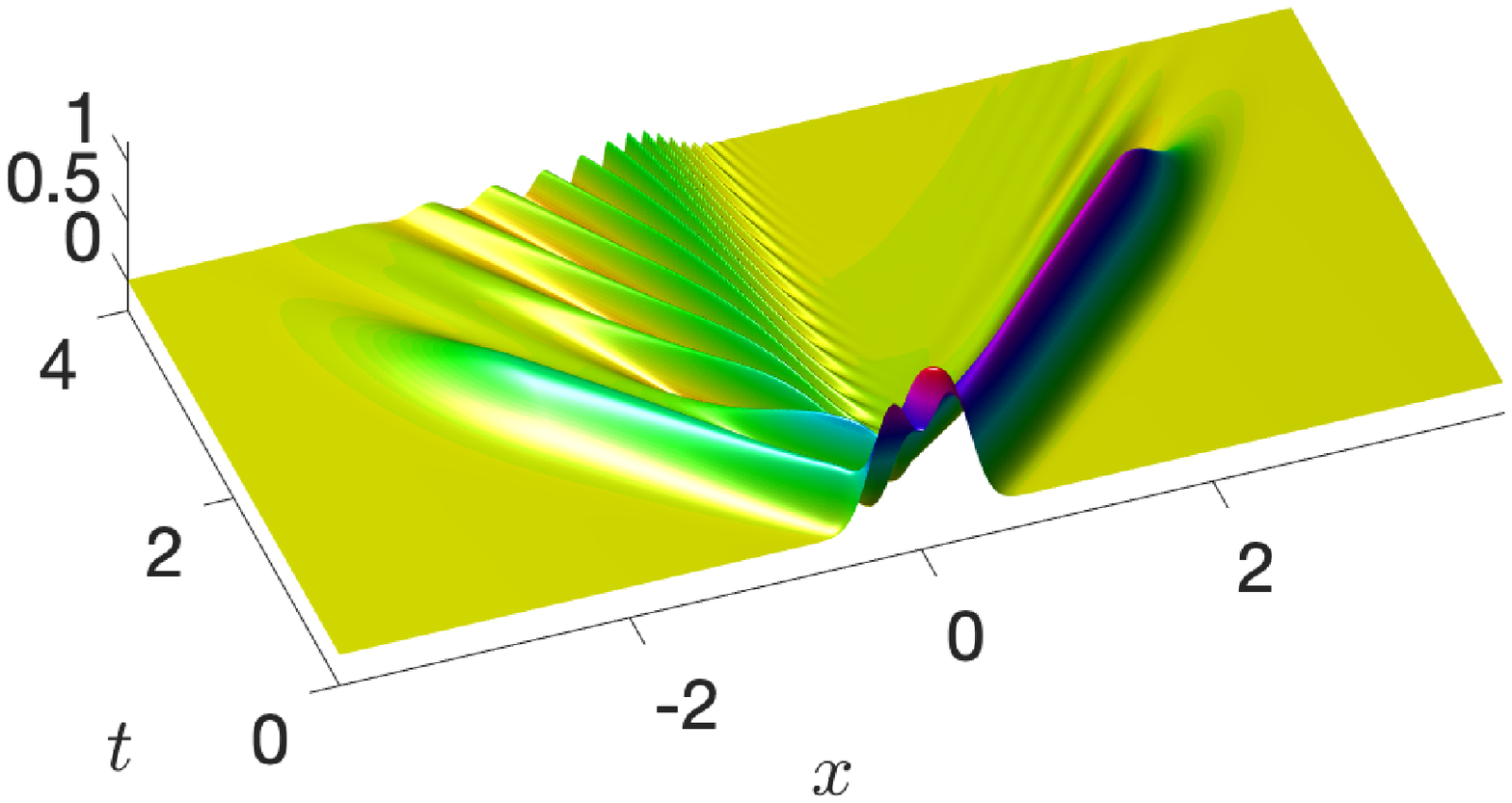}
 \caption{numerical solution $\delta=0.2$}
 \label{fig:ex5p3solutionsb}
\end{subfigure}
\begin{subfigure}{.32\textwidth}
 \centering
	\includegraphics[width=1\textwidth]{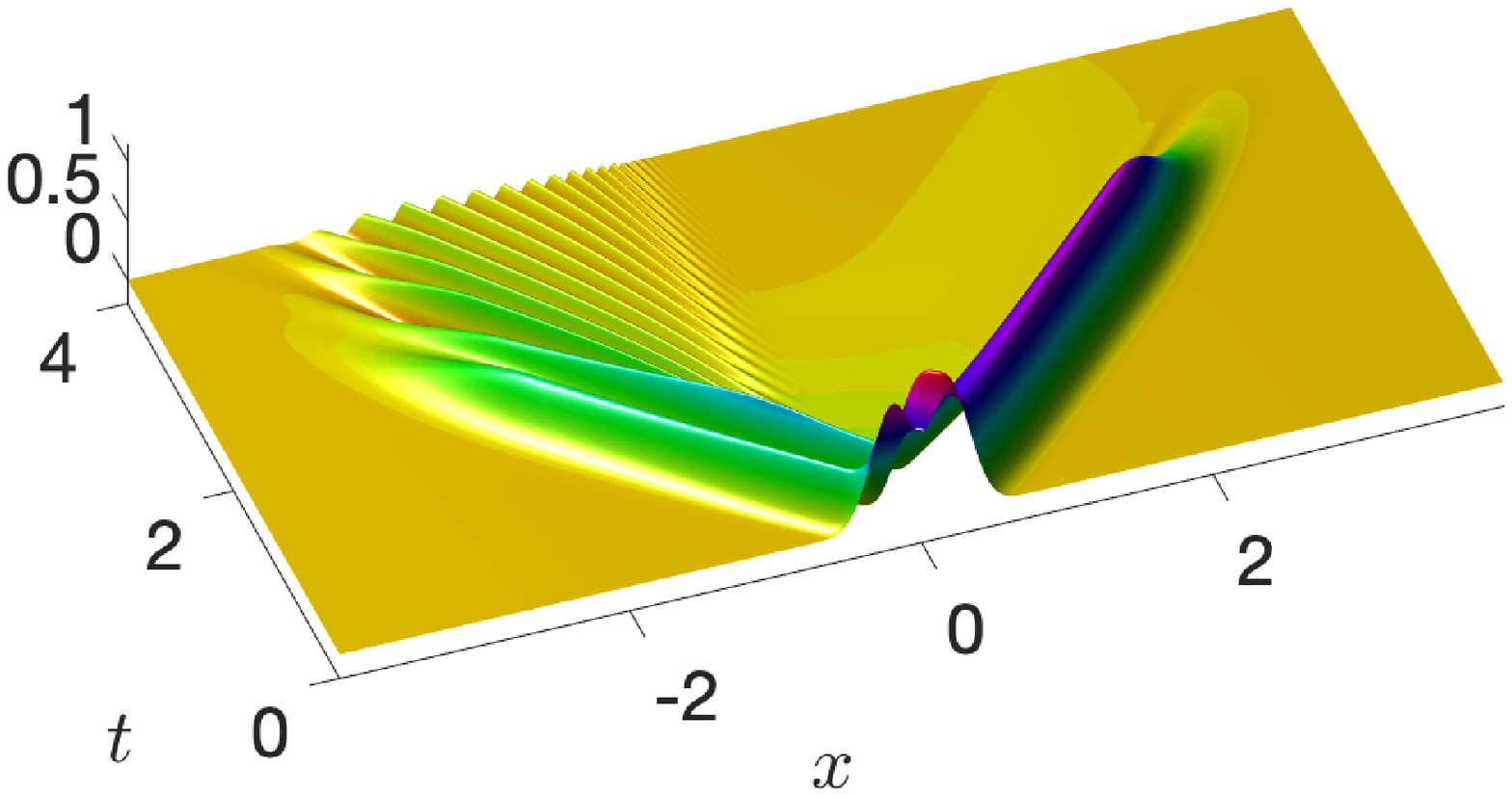}
 \caption{numerical solution $\delta=0.1$}
 \label{fig:ex5p3solutionsc}
 \end{subfigure}
 
 \begin{subfigure}{.32\textwidth}
 \centering
	\includegraphics[width=1\textwidth]{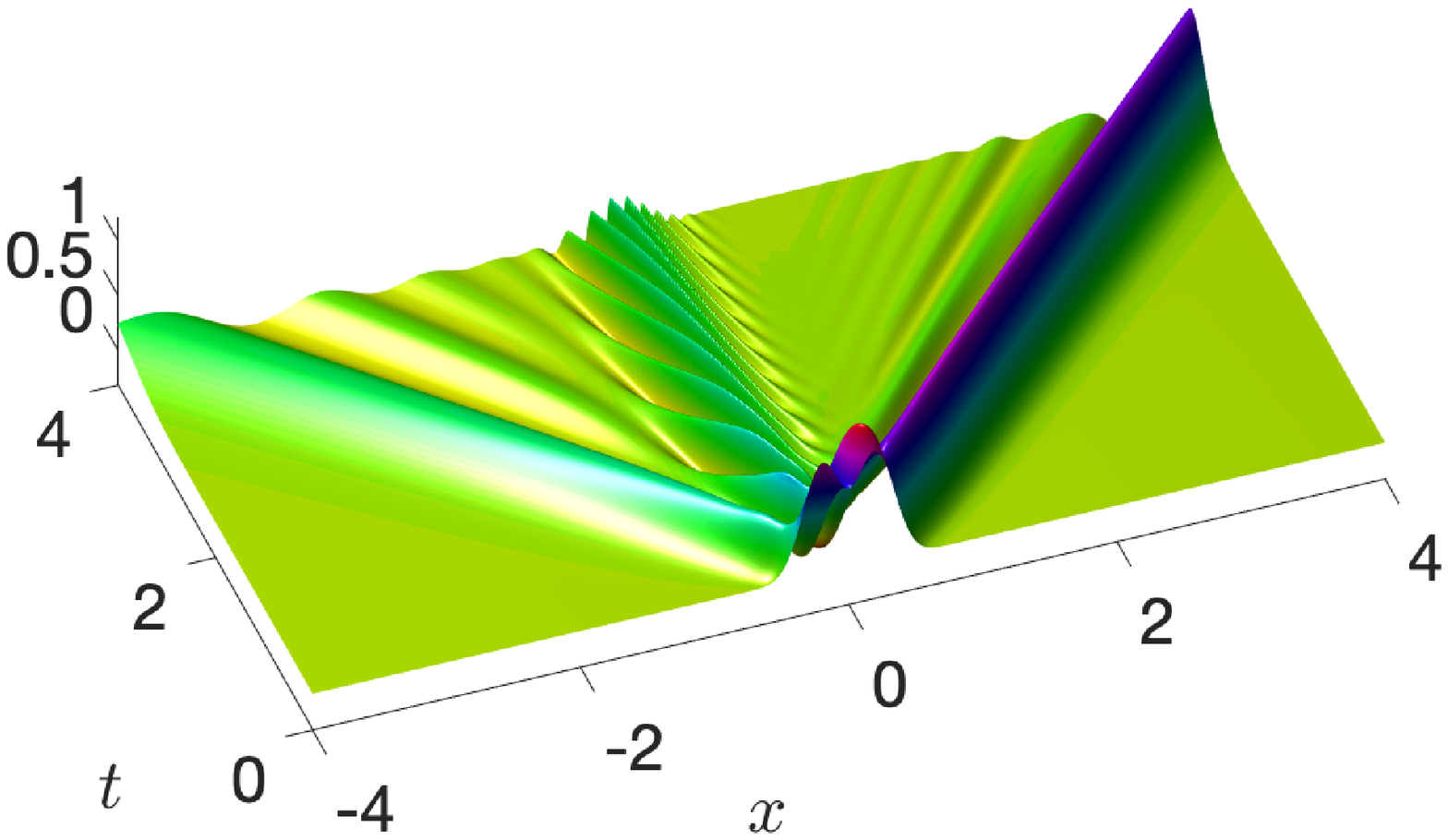}
 \caption{reference solution $\delta=0.3$}
 \label{fig:ex5p3solutionsd}
\end{subfigure}
\begin{subfigure}{.32\textwidth}
 \centering
	\includegraphics[width=1\textwidth]{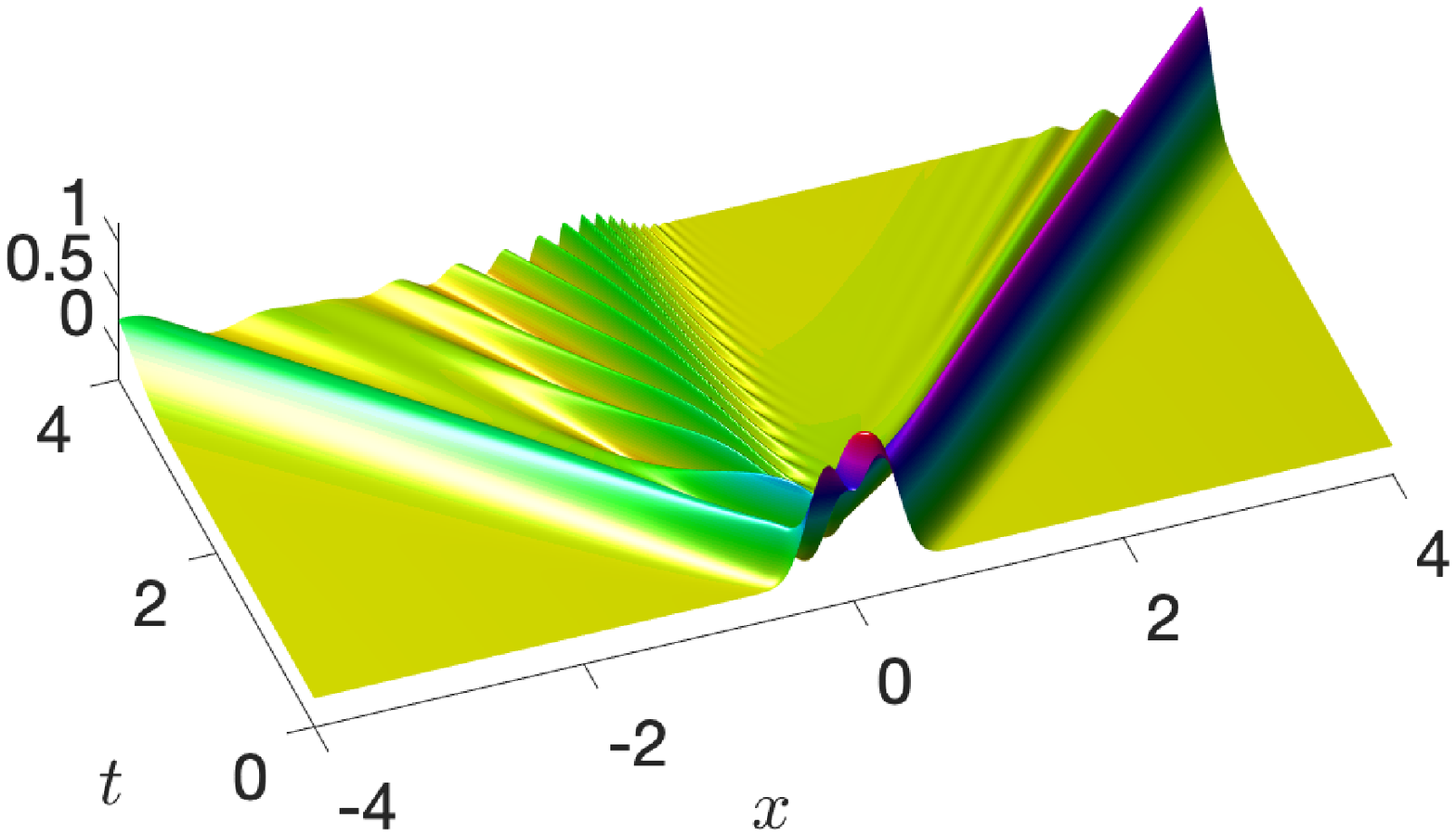}
 \caption{reference solution $\delta=0.2$}
 \label{fig:ex5p3solutionse}
\end{subfigure}
\begin{subfigure}{.32\textwidth}
 \centering
	\includegraphics[width=1\textwidth]{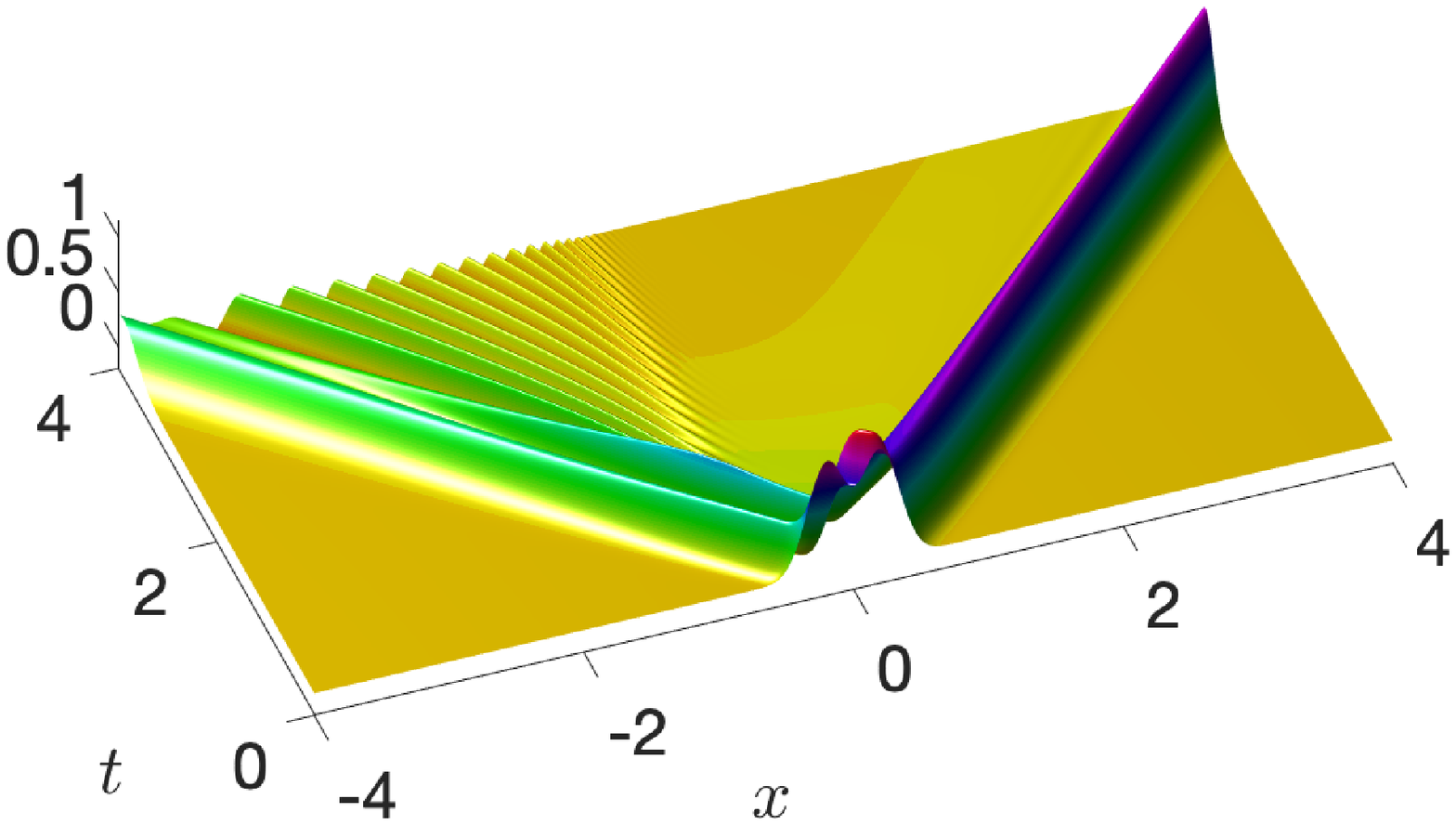}
 \caption{reference solution $\delta=0.1$}
 \label{fig:ex5p3solutionsf}
 \end{subfigure}
	\caption{(Example 5.3) Numerical solutions for $\delta=0.3,0.2,0.1$ and reference solutions up to $T=4$. The numerical solutions are obtained by taking $h=2^{-8}$ and $m=400$.}
	\label{fig:ex5p3solutions}
\end{figure}

\begin{figure}[htbp]
\centering
\begin{subfigure}{.32\textwidth}
 \centering
	\includegraphics[width=1\textwidth]{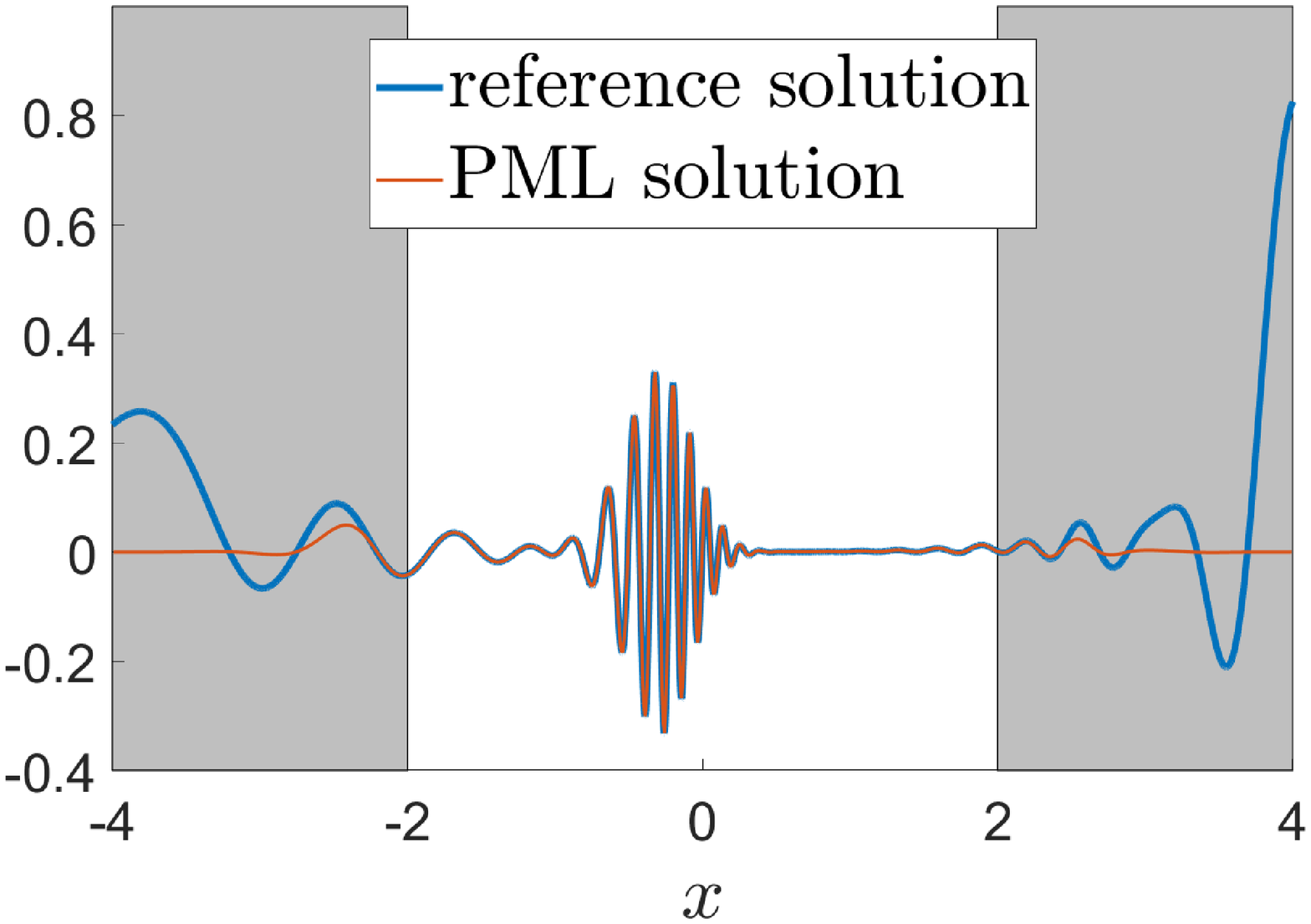}
 \caption{$\delta=0.3$}
 \label{fig:ex5p3solutionsT2a}
\end{subfigure}
\begin{subfigure}{.32\textwidth}
 \centering
	\includegraphics[width=1\textwidth]{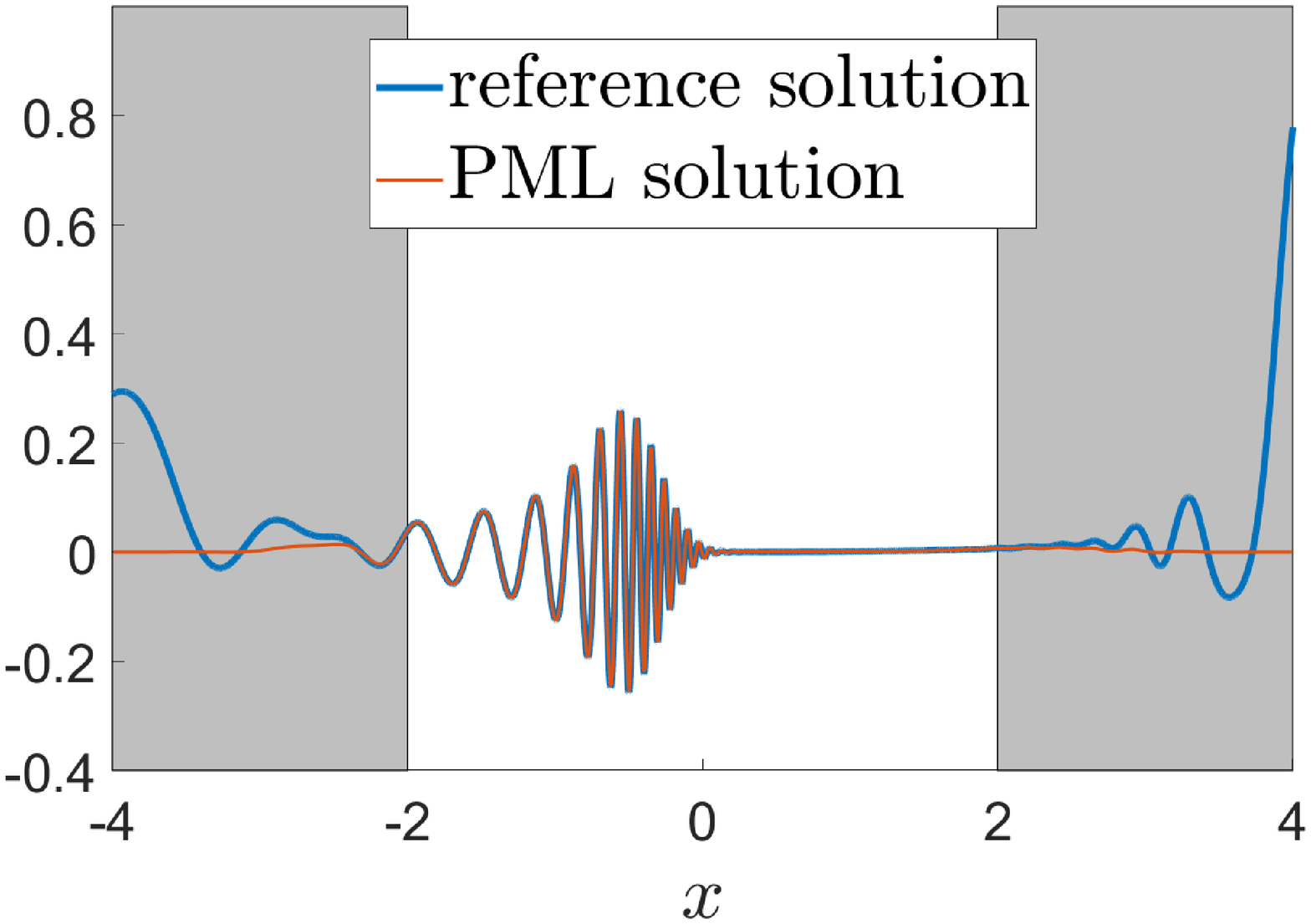}
 \caption{$\delta=0.2$}
 \label{fig:ex5p3solutionsT2b}
\end{subfigure}
\begin{subfigure}{.32\textwidth}
 \centering
	\includegraphics[width=1\textwidth]{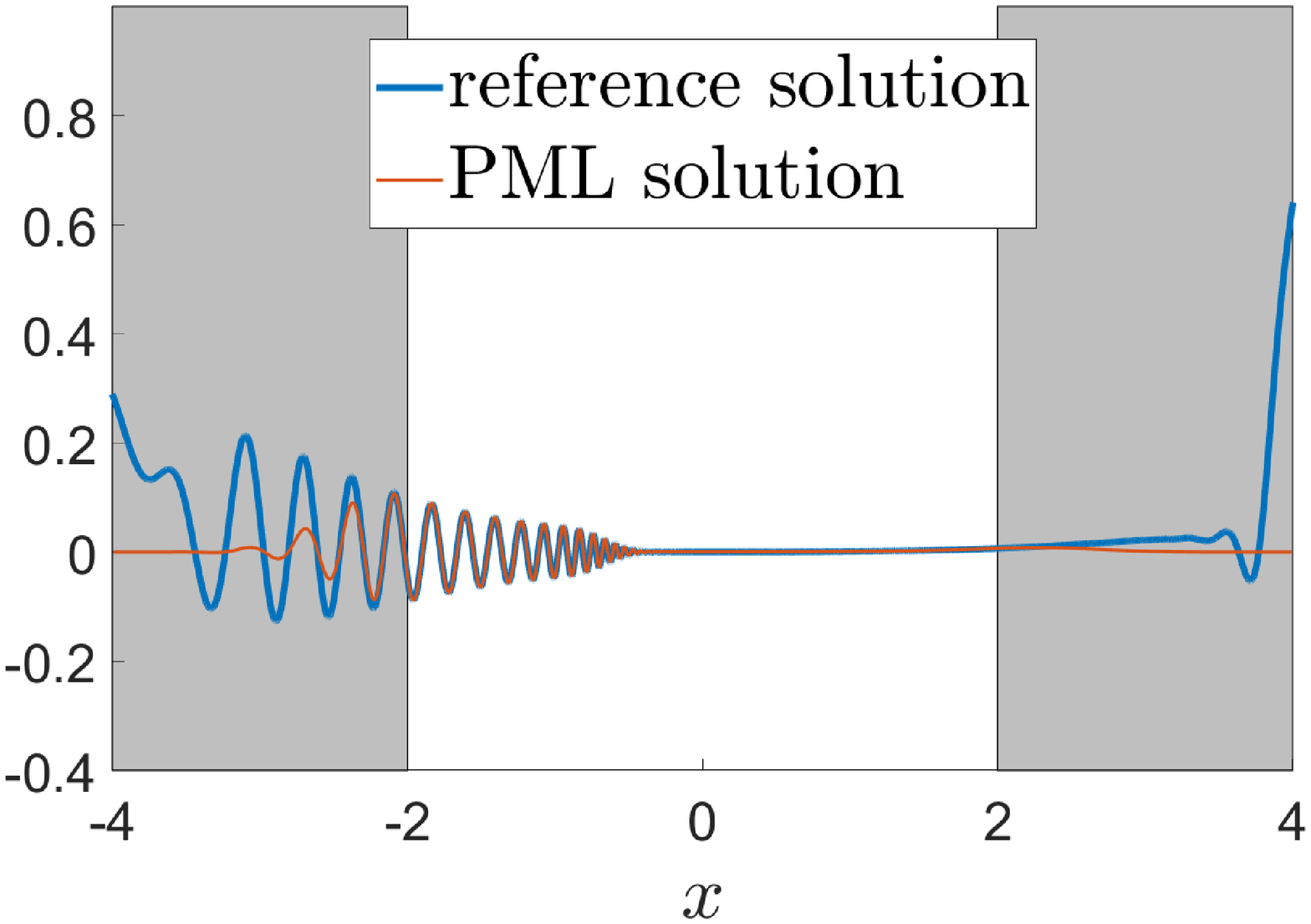}
 \caption{$\delta=0.1$}
 \label{fig:ex5p3solutionsT2c}
 \end{subfigure}
	\caption{(Example 5.3) The comparison of the numerical solution and reference solution at time $t=4$. The PML layers for numerical solutions are shaded in light grey.}
	\label{fig_ex5p3solutionsT4}
\end{figure}

Figure~\ref{fig:ex5p3solutions} shows the evolution of numerical and reference solutions by taking $T=4$, $h=2^{-8}$, $\tau=1/12000$ and $m=400$ for different $\delta=0.3,0.2,0.1$. Figure~\ref{fig_ex5p3solutionsT4} shows numerical and reference solutions at $t=4$, which again indicates the numerical solutions decay exponentially in PML layers.

Table~\ref{tab_5p3_l2} lists the errors $e_h$ at $t=4$ and the spatial convergence order for numerical solutions, which verify the second-order accuracy of the spatial discretization by refining $h$. Table~\ref{tab_5p3_dh_l2} lists the errors $e_\delta$ and $\delta$-convergence orders between the numerical solutions and exact solutions of local problem \eqref{eq_localinPML1}--\eqref{eq_localinPML1} by vanishing $\delta$ and $h$ simultaneously at $t=2$. The behavior of the errors is similar to that shown in the examples above for sufficiently small $\delta$.

\begin{table}
\centering
\begin{tabular}{c|cc|cc|cc}
\hline
$h$	&	$\delta=0.3$  & order & $\delta=0.2$ & order &  $\delta=0.1$ & order\\ \hline
$2^{-5}$	& 3.74e-02 &  --& 4.7483e-02 & -- &  4.35e-02  & --\\
$2^{-6}$ &  9.96e-03  & 1.91 &  1.3050e-02  & 1.86 &  1.19e-02 & 1.87\\
$2^{-7}$ & 2.52e-03 & 1.98 &  3.3125e-03& 1.98 &  2.99e-03 & 2.00\\
$2^{-8}$ & 6.36e-04 & 1.98 & 8.25e-04  & 2.01 & 7.21e-04  & 2.05 \\\hline
\end{tabular}
\caption{(Example 5.3) errors $e_h$ and convergence orders at $t=4$.}\label{tab_5p3_l2}
\end{table}

\begin{table}
\centering
\begin{tabular}{c|cc|cc|cc}
\hline
$h$	&	$\delta=h$  & order & $\delta=2h$ & order &  $\delta=3h$ & order\\ \hline
$2^{-5}$	&  5.3521e-02 &  -- & 1.1867e-01 & -- & 1.4328e-01 & --\\
$2^{-6}$ &  1.3937e-02& 1.94 &  4.2328e-02 & 1.49 & 8.1629e-02 & 0.81\\
$2^{-7}$ &  3.4432e-03  & 2.02 &  1.0714e-02 & 1.98 &  2.2952e-02 & 1.83\\
$2^{-8}$ &  8.3835e-04  & 2.03 &  2.6479e-03& 2.01&  5.6996e-03 & 2.01\\\hline
\end{tabular}
\caption{(Example 5.3) The errors $e_\delta$ and $\delta$-convergence orders between numerical solutions and exact solutions of local PML problem \eqref{eq_localinPML1}--\eqref{eq_localinPML2} by vanishing $\delta$ and $h$ simultaneously at $t=2$.}\label{tab_5p3_dh_l2}
\end{table}

\section{Conclusion}
The nonlocal PML and its numerical discretization of a nonlocal wave equation in unbounded spatial domains are studied in this paper. We first propose the nonlocal PML equation with a time-dependent nonlocal operator which has a complex-valued kernel and is defined by a convolution. This feature differs from the results in \cite{wildman2011,wildman2012a}. After discretizing the contour integrals arising from the inverse Laplace transform, applying the AC scheme and introducing auxiliary functions, we get a truncated semidiscrete nonlocal PML problem which only involves a finite number of degrees of freedom. A Verlet-type ODE solver is used for the time integration. Numerical experiments demonstrate the effectiveness and accuracy of our proposed PML.


\section*{Appendix A.} \label{sec_Ap}
\textbf{We first consider $\Omega_{\mathcal{K}}$ for the kernel~\eqref{eq_ker1}.} First, we give the analytic continuation of the reference kernel $\gamma_0(s) = \frac{1}{2c_0^3} e^{-\frac{|s|}{c_0}}$ by
$$\gamma_0(\tilde x-\tilde y) = \frac{1}{2c_0^3} e^{-\frac{\rho(\tilde x,\tilde y)}{c_0}},$$
where $\rho(\tilde x,\tilde y)=\sqrt{(\tilde x-\tilde y)^2}$ is its analytic branch $\Re \rho(\tilde x,\tilde y)\geq0$. Let $z=z_1+\i z_2$ and $s = s_1 + \i s_2$ where $z_1,z_2,s_1,s_2\in\R$. We have
\begin{align*}
\tilde x-\tilde y &= (x-y) + \frac{z}{s}\Big(\int_0^x\sigma(t)\mathrm dt-\int_0^y\sigma(t)\mathrm dt\Big)\\
&= (x-y) \Big[\Big( 1 + \frac{z_1s_1+z_2s_2}{|s|^2}\cdot \frac{\int_0^x\sigma(t)\mathrm dt-\int_0^y\sigma(t)\mathrm dt}{x-y}\Big) x+\i\frac{-z_1s_2+z_2s_1}{|s|^2} \cdot\frac{\int_0^x\sigma(t)\mathrm dt-\int_0^y\sigma(t)\mathrm dt}{x-y}\Big].
\end{align*}
Set
\begin{equation}\label{Ge}
g=\frac{\int_0^x\sigma(t)\mathrm dt-\int_0^y\sigma(t)\mathrm dt}{x-y}.
\end{equation}
 It is clear that $0\leq g\leq1$ since $\sigma\leq 1$. Therefore, the real part of $\tilde x-\tilde y$ is
\begin{align*}
 (x-y)\Big(1 + \frac{(z_1s_1+z_2s_2)g}{|s|^2} \Big) &=  (x-y)\Big(\frac{|s|^2 + (z_1s_1+z_2s_2)g}{|s|^2} \Big)\\
 &= \frac{x-y}{|s|^2}\Big[ \Big( s_1+\frac{z_1g}{2} \Big)^2 +  \Big( s_2+\frac{z_2g}{2} \Big)^2 - \frac14|z|^2g^2\Big],
\end{align*}
which implies that $\rho(\tilde x,\tilde y)=\sqrt{(\tilde x-\tilde y)^2}$ as a function of $s$ is analytic in the $s$-domain
\begin{align*}
\Omega_{\mathcal{K}}:\ \Big( s_1+\frac{z_1}{2} \Big)^2 +  \Big( s_2+\frac{z_2}{2} \Big)^2 > \frac14|z|^2, \quad \forall x,y. 
\end{align*}

\noindent\textbf{We then consider Talbot's contour parameters $\mu$ and $\nu$ for the Gaussian kernel~\eqref{eq_gaukernel}.} The analytic continuation of the kernel $\gamma(y-x,\frac{x+y}{2}) =  \frac{4}{\delta^3}\sqrt{\frac{10^3}{\pi}} e^{-10\frac{(x-y)^2}{\delta^2}}$ is given by
\begin{align}
\gamma(\tilde y-\tilde x,\frac{\tilde x+\tilde y}{2}) =  \frac{4}{\delta^3}\sqrt{\frac{10^3}{\pi}} e^{-10\frac{(\tilde x-\tilde y)^2}{\delta^2}}.
\end{align}
Note the $\Omega_\mathcal{K}$ is the whole complex plane for any given $z\in\mathbb{C}$ and $x,y\in\R$. However, to ensure the stability, we have to choose the Talbot's contour parameters $\mu$ and $\nu$ such that $\Re[(\tilde x-\tilde y)^2] \geq 0$. Denote by $\zeta=\zeta_1+\i \zeta_2=\frac{z}{\xi_j}$. We have
\begin{align*}
\Re\Big[(\tilde x-\tilde y)^2\Big] =& \Re\Big[\Big( (x+\zeta_1\int_0^x\sigma(t)\mathrm dt+\i\zeta_2\int_0^x\sigma(t)\mathrm dt) - (y+\zeta_1\int_0^y\sigma(t)\mathrm dt+\i\zeta_2\int_0^y\sigma(t)\mathrm dt) \Big)^2\Big] \\
=& \Big( (x+\zeta_1\int_0^x\sigma(t)\mathrm dt)- (y+\zeta_1\int_0^y\sigma(t)\mathrm dt)\Big)^2 -\Big (\zeta_2\int_0^x\sigma(t)\mathrm dt-\zeta_2\int_0^y\sigma(t)\mathrm dt\Big)^2\\
=&(x-y)^2\big[ (1+\zeta_1 g)^2 - \zeta_2^2 g^2\big]\\
=&(x-y)^2\big[ \big(1+(\zeta_1-\zeta_2)g\big) \big(1+(\zeta_1+\zeta_2)g\big)\big],
\end{align*}
where $g$ is defined in \eqref{Ge} with $0\leq g\leq1$ as $0\leq\sigma\leq 1$. 

To ensure $\Re[(\tilde x-\tilde y)^2] \geq 0$, we have $\big(1+(\zeta_1-\zeta_2)g\big) \big(1+(\zeta_1+\zeta_2)g\big)\geq0$ for any $g\in[0,1]$, which implies that $\zeta_1-\zeta_2\geq -1$ and $\zeta_1+\zeta_2\geq-1$. Therefore, we may simply choose Talbot's contour parameters $\mu$ and $\nu$ such that
\begin{align}
\sqrt{2}|z|\leq |\xi_j|,\quad \forall\ j=1,2,\cdots,m.
\end{align}


\bibliographystyle{siam}
\bibliography{../referrence}

\end{document}